\documentclass[11pt,a4paper]{article}
\usepackage[T1]{fontenc}
\usepackage{lmodern}
\usepackage[margin=25mm]{geometry}
\usepackage{amsmath,amssymb,amsthm}
\usepackage{enumitem}
\usepackage{tikz}
\usepackage[colorlinks=true,linkcolor=blue,citecolor=red,urlcolor=blue]{hyperref}

\newcommand{\Z}{\mathbb{Z}}
\newcommand{\F}{\mathbb{F}}
\newcommand{\C}{\mathbb{C}}
\newcommand{\ket}[1]{\lvert #1\rangle}

\theoremstyle{plain}
\newtheorem{theorem}{Theorem}[section]
\newtheorem{lemma}[theorem]{Lemma}
\newtheorem{proposition}[theorem]{Proposition}
\newtheorem{corollary}[theorem]{Corollary}
\theoremstyle{definition}
\newtheorem{definition}[theorem]{Definition}
\newtheorem{example}[theorem]{Example}
\newtheorem{question}[theorem]{Question}

\setlist[enumerate,1]{label=\textup{(\roman*)},leftmargin=*}
\title{\bf Extremal permutations of $\mathbb{Z}_n$ with many good slopes
and applications to quantum Latin squares\thanks{This work has been granted by the Juan de la Cierva 2024 grant (JDC2024-053082-I), funded by Ministerio de Ciencia, Innovación y Universidades (MICIU/AEI/10.13039/501100011033) and co-funded by the European Social Fund Plus (FSE+).}}
\author{
Dipak K. Bhunia\\
\small Departament de Matem\`atiques,\\ \small Universitat Polit\`ecnica de Catalunya,\\
\small Barcelona, Spain
}
\date{}

\begin{document}
\maketitle

\begin{abstract}
For $n\ge2$, a good slope of a permutation $f$ of $\mathbb{Z}_n=\mathbb{Z}/n\mathbb{Z}$ is an element $m\in\mathbb{Z}_n$ for which $x\mapsto f(x)+mx$ is a permutation. We prove the sharp upper bound $\varphi(n)$ on the number of good slopes, where $\varphi$ is Euler's totient function. We classify and count the equality cases for prime powers and products of two distinct odd primes. Non-affine equality cases exist for every composite $n\ge6$. We also prove the sharp bound $p_0-1$ for the size of a set of good slopes with invertible differences between distinct elements, where $p_0$ is the least prime divisor of $n$. Finally, we establish gap results below $\varphi(n)$ and apply our bounds to a construction of quantum Latin squares due to Ball and Simoens.

\medskip
\noindent\textbf{Keywords.} Good slopes; permutations of residue rings;
extremal permutations; congruence preservation; mutually orthogonal quantum Latin squares.

\medskip
\noindent\textbf{MSC 2020.} 05A05, 05B15, 81P70.
\end{abstract}

\section{Introduction}

Let $n\ge2$, write $\Z_n=\Z/n\Z$, and let $\Z_n^\times$ be
its group of units, that is, its invertible elements. For a permutation
$f:\Z_n\to\Z_n$, let $D_n(f)$ be the set of all $m\in\Z_n$ such that
$x\mapsto f(x)+mx$ is a permutation. We call these residues
\emph{good slopes}. Since $f$ is a permutation, $0\in D_n(f)$. A
subset of $D_n(f)$ is \emph{admissible} if the difference of any two
distinct elements is a unit. We write $\mu_n(f)$ for the largest size
of an admissible subset. Let $\varphi$ be Euler's totient function and
$p_0(n)$ the least prime divisor of $n$. An affine permutation has the
form $f(x)=ax+b$, where $a\in\Z_n^\times$ and $b\in\Z_n$. Its set of
good slopes is $D_n(f)=\Z_n^\times-a=\{u-a:\ u\in\Z_n^\times\}$.
Hence, it has exactly $\varphi(n)$ good slopes.

The problem connects permutation theory with finite geometry.
Good slopes are defined in the same way for a map $h$ from a finite
field to itself. A slope $a$ is \emph{determined} by the graph of $h$
if $h(x)-h(y)=a(x-y)$ for some $x\ne y$. Thus, $a$ is not determined
precisely when $x\mapsto h(x)-ax$ is a permutation, that is, when $-a$
is a good slope of $h$. The direction theorems of R\'edei, and later of
Blokhuis, Ball, Brouwer, Storme and Sz\H{o}nyi, give strong
restrictions on maps with many good slopes \cite{Ball,BBBSS,BBS,Redei}.
De Beule surveys direction problems in affine spaces \cite{DeBeule}.

There is also a connection with complete mappings and
orthomorphisms. These are permutations $f$ for which $f(x)+x$ and
$f(x)-x$, respectively, are permutations \cite{Evans}. In our notation,
they are the permutations with $1\in D_n(f)$ and with
$-1\in D_n(f)$, respectively. Thus, these classical notions concern a
single slope, whereas we study the full set of good slopes of a
permutation.

Congruence-preserving maps play a central role in this paper. A map
$g:\Z_n\to\Z_n$ is \emph{congruence-preserving} if $x\equiv y\pmod d$
implies $g(x)\equiv g(y)\pmod d$ for every divisor $d$ of $n$.
Stones and Wanless \cite{SW} call orthomorphisms with this property
compatible. They also study $d$-compound orthomorphisms, which
preserve congruences modulo a single divisor $d$ of $n$.
C\'egielski, Grigorieff and Guessarian characterise the
congruence-preserving maps between residue rings by means of
polynomials with rational coefficients \cite{CGG}.

In this paper, we prove that every permutation satisfies
$|D_n(f)|\le\varphi(n)$ and $\mu_n(f)\le p_0(n)-1$.
We call a permutation \emph{extremal} when the first bound is attained.
For every composite $n\ge6$, some non-affine permutation attains both
bounds simultaneously. We classify and count all extremal permutations
when $n$ is a prime power or a product of two distinct odd primes.
For these moduli, every extremal permutation has the form
$f(x)=ax+b+g(x)$, where $a$ is a unit, $b\in\Z_n$ and $g$ is
\emph{strictly triangular} with respect to a chain of ideals
$\Z_n=I_0\supsetneq I_1\supsetneq\cdots\supsetneq I_r=\{0\}$. This
means that $x-y\in I_k$ implies $g(x)-g(y)\in I_{k+1}$ for $0\le k<r$.
For every other modulus with at least two distinct prime divisors, we
construct an extremal permutation outside this class.

Let $f$ be a congruence-preserving permutation of $\Z_{p^e}$, where
$p$ is a prime, and write $x=\sum_{j<e}x_jp^j$ and
$f(x)=\sum_{j<e}y_jp^j$ in base $p$. Then, $y_j$ depends only on
$x_0,\dots,x_j$. Once $x_0,\dots,x_{j-1}$ are fixed, $y_j$ is obtained
from $x_j$ by a permutation of the set $\{0,\dots,p-1\}$ of digits. We
describe good slopes through these \emph{digit maps}. This gives a
classification of the congruence-preserving permutations of $\Z_n$
with more than $n/2$ good slopes. For permutations induced by
polynomials, we make this classification explicit.

We also study the gap below the maximum. If $n$ is odd and composite
and $q$ is its largest prime divisor, put $G(n)=\varphi(n)(q-2)/(q-1)$.
Every non-extremal congruence-preserving permutation has at most
$G(n)$ good slopes. The same bound holds for every non-extremal
permutation when $n$ is a product of two distinct odd primes, and when
$n=mq$ for a prime $q>m$ with $\varphi(m)\le(q-1)/2$. A reduction
through preserved congruences gives further cases whenever the smaller
moduli are settled.
We construct permutations attaining $G(n)$ for every odd composite $n$.
In a companion paper \cite{Gaps}, we study the permutations of
prime-power order that are not congruence-preserving, and we prove the
gap bound for further moduli.

The final application concerns quantum Latin squares. A quantum Latin
square of order $n$ is an $n\times n$ array of unit vectors of $\C^n$
whose rows and columns are orthonormal bases. These squares were
introduced by Musto and Vicary \cite{MV2016}. Two such squares are
orthogonal if the tensor products of their corresponding entries form
an orthonormal basis of $\C^n\otimes\C^n$. Goyeneche, Raissi,
Di Martino and \.Zyczkowski proposed a notion of orthogonality for
quantum Latin squares \cite{GRMZ}, and Musto and Vicary proved that it
is equivalent to this tensor-product condition \cite{MV2019}.
The construction of Ball and Simoens \cite{BS} produces pairwise
orthogonal quantum Latin squares from a non-affine permutation and an
admissible set of good slopes. Our results determine the largest
families supplied by this construction over $\Z_n$.

The paper is organised as follows. Section~\ref{sec:prelim} collects
the notation and preliminary results. Section~\ref{sec:universal}
proves the universal bounds, and Section~\ref{sec:constructions}
determines their sharpness for non-affine permutations.
Sections~\ref{sec:primepowers} and~\ref{sec:mixed} treat the
classification, enumeration and construction of extremal permutations.
Section~\ref{sec:half} studies congruence-preserving permutations,
including polynomial permutations. Section~\ref{sec:gap} establishes
gap results below the maximum. Section~\ref{sec:redei} discusses the
R\'edei polynomial over $\Z_n$.
Section~\ref{sec:qls} gives the application to quantum Latin squares.
Finally, Section~\ref{sec:questions} presents further questions.

\section{Preliminaries}\label{sec:prelim}

In this section, we fix our notation and collect basic facts on good
slopes and admissible sets. We also show how good slopes behave under
a preserved congruence, and we recall two elementary facts on
polynomials and power sums, which are used in Section~\ref{sec:universal}.

Throughout, $n\ge2$ is an integer, and $\Z_n^\times$, $\varphi$ and
$p_0(n)$ are as in the introduction. For a prime $p$ and an integer
$e\ge1$, we write $p^e\,\|\,n$ if $p^e\mid n$ and $p^{e+1}\nmid n$. By the
Chinese remainder theorem, we identify $\Z_n$ with
$\prod_{p^e\|n}\Z_{p^e}$ whenever this is convenient.
For a statement $P$, we put
$$
[P]=\begin{cases}1&\text{if $P$ holds},\\0&\text{otherwise}.\end{cases}
$$
For a subset $S$ of a ring and elements $a$ and $u$ of this ring, we
write
$$
S+a=\{s+a:\ s\in S\},\quad uS=\{us:\ s\in S\}\quad\text{and}\quad
{-S}=\{-s:\ s\in S\},
$$
and we put $S-a=S+(-a)$.
For instance, $a+p\Z_{p^e}$ is the residue class of $a$ modulo $p$
in $\Z_{p^e}$.
Whenever we write $\Z_{p^e}$, the integer $p$ is a prime and $e\ge1$.
Further restrictions, such as $e\ge2$ or $p$ odd, are stated where
they are needed. For a prime power $q$, we write $\F_q$ for the field
with $q$ elements. We identify $\Z_p$ with $\F_p$. For $e\ge2$, the
rings $\Z_{p^e}$ and $\F_{p^e}$ are not isomorphic.

A map $g:\Z_n\to\Z_n$ is \emph{affine} if $g(x)=ax+b$ for some
$a,b\in\Z_n$. The additive maps of the cyclic group $\Z_n$ are the maps
$x\mapsto ax$. Thus, a map of $\Z_n$ is affine exactly when it is the sum
of an additive map and a constant. An affine map $x\mapsto ax+b$ is a
permutation if and only if $a\in\Z_n^\times$. A map is \emph{polynomial} if
it is induced by a polynomial in $\Z[X]$. Let $d$ be a divisor of $n$. A
map $g$ of $\Z_n$ \emph{preserves congruences modulo $d$} if
$x\equiv y\pmod d$ implies $g(x)\equiv g(y)\pmod d$. It is
\emph{congruence-preserving} if it preserves congruences modulo every
divisor of $n$. Polynomial maps are congruence-preserving, because $x-y$
divides $F(x)-F(y)$ for every $F\in\Z[X]$.

\begin{definition}\label{def:slopes}
Let $f$ be a permutation of $\Z_n$. An element $m\in\Z_n$ is a \emph{good
slope} of $f$ if the map $x\mapsto f(x)+mx$ is a permutation of $\Z_n$. We
write $D_n(f)$ for the set of good slopes of $f$. A subset
$M\subseteq D_n(f)$ is \emph{admissible} if $m-m'\in\Z_n^\times$ for all
distinct $m,m'\in M$. We write $\mu_n(f)$ for the largest size of an
admissible subset of $D_n(f)$.
\end{definition}

Since $f$ is a permutation, $0\in D_n(f)$ and $\mu_n(f)\ge1$.
Geometrically, $m\in D_n(f)$ means that each of the $n$ lines
$\{(x,y):y+mx=c\}$ with $c\in\Z_n$ contains exactly one point of the graph
of $f$. When $n$ is prime, every subset of $D_n(f)$ is admissible and
$\mu_n(f)=|D_n(f)|$. If $f(x)=ax+b$ is an affine permutation, then
$f(x)+mx=(a+m)x+b$. Hence, $D_n(f)=\Z_n^\times-a$ and
$|D_n(f)|=\varphi(n)$.

\begin{lemma}\label{lem:normal}
Let $f$ be a permutation of $\Z_n$, let $a\in D_n(f)$ and $b,c\in\Z_n$, and
let $u,v\in\Z_n^\times$. Then, $x\mapsto f(x)+ax+b$ and
$x\mapsto uf(vx+c)$ are permutations of $\Z_n$, and
$$
D_n\bigl(f(x)+ax+b\bigr)=D_n(f)-a,\quad
D_n\bigl(uf(vx+c)\bigr)=uv\,D_n(f).
$$
Both operations preserve the number of good slopes, the value of $\mu_n$,
and the property of being affine.
\end{lemma}

\begin{proof}
For $m\in\Z_n$, we have $f(x)+ax+b+mx=f(x)+(a+m)x+b$, and adding the
constant $b$ does not affect bijectivity. Hence, $m$ is a good slope of
$f(x)+ax+b$ if and only if $a+m\in D_n(f)$. Taking $m=0$, we see that
$f(x)+ax+b$ is a permutation. Now put $g(x)=uf(vx+c)$ and $y=vx+c$. The map
$x\mapsto y$ is a permutation of $\Z_n$, and
$g(x)+mx=u\bigl(f(y)+u^{-1}v^{-1}my\bigr)-v^{-1}mc$. Therefore,
$x\mapsto g(x)+mx$ is a permutation if and only if
$u^{-1}v^{-1}m\in D_n(f)$. Taking $m=0$ again shows that $g$ is a
permutation. Translating a set, or multiplying it by a unit, preserves its
size and the invertibility of its differences. Finally, both operations
can be reversed by operations of the same kind, and they send affine maps
to affine maps.
\end{proof}

The next lemma describes admissible sets through the residues of their
elements modulo the prime divisors of $n$.

\begin{lemma}\label{lem:admissible}
\begin{enumerate}
\item Let $m,m'\in\Z_n$. Then, $m-m'\in\Z_n^\times$ if and only if
$m\not\equiv m'\pmod p$ for every prime divisor $p$ of $n$.
\item For every prime divisor $p$ of $n$, let $T_p$ be a nonempty
subset of the set of residues modulo $p$, and let $D$ be the set of $m\in\Z_n$ with
$m\bmod p\in T_p$ for every prime $p\mid n$. Then,
$|D|=n\prod_{p\mid n}(|T_p|/p)$. Moreover, the largest subsets of $D$ whose
pairwise differences are units have exactly $\min_{p\mid n}|T_p|$
elements.
\item Let $p$ be prime, let $e\ge1$, and let $f$ be a permutation
of $\Z_{p^e}$. Then,
$$
\mu_{p^e}(f)=\bigl|\{m\bmod p:m\in D_{p^e}(f)\}\bigr|.
$$
\end{enumerate}
\end{lemma}

\begin{proof}
An element of $\Z_n$ is a unit if and only if it is divisible by no prime
divisor of $n$. This proves (i). For (ii), let $r$ be the product of the
distinct prime divisors of $n$. By the Chinese remainder theorem, reduction
modulo the primes $p\mid n$ is a surjective map from $\Z_n$ onto
$\prod_{p\mid n}\Z_p$, and each of its fibres has $n/r$ elements. The set
$D$ is the preimage of $\prod_{p\mid n}T_p$. Hence,
$|D|=(n/r)\prod_{p\mid n}|T_p|$, which is the stated value. By (i), the
elements of a subset of $D$ with pairwise invertible differences have
distinct residues modulo every $p\mid n$, and these residues lie in $T_p$.
Thus, such a subset has at most $\min_p|T_p|$ elements. Conversely, put
$t=\min_p|T_p|$, and choose distinct residues $c_{p,1},\dots,c_{p,t}\in T_p$
for every $p\mid n$. The Chinese remainder theorem gives
$m_1,\dots,m_t\in\Z_n$ with $m_k\equiv c_{p,k}\pmod p$ for all $p$ and $k$.
These elements lie in $D$, and by (i) their pairwise differences are
units.

For (iii), the elements of an admissible subset must have distinct
residues modulo $p$, by (i). Conversely, choose one good slope from
each residue class modulo $p$ that meets $D_{p^e}(f)$. The chosen
slopes have pairwise unit differences, again by (i). This proves the
formula.
\end{proof}

Next, we describe the good slopes of a permutation that preserves
congruences modulo a divisor $d$ of $n$. Orthomorphisms with this
property are the $d$-compound orthomorphisms of \cite{SW}.

\begin{lemma}\label{lem:split-general}
Let $n=dk$, where $d,k\ge2$, and let $f$ be a permutation of $\Z_n$
that preserves congruences modulo $d$. Define $g:\Z_d\to\Z_d$ by
$g(x\bmod d)=f(x)\bmod d$. For $0\le r<d$, define
$f_r:\Z_k\to\Z_k$ by
$$
f(r+dt)-f(r)=d f_r(t),
$$
where multiplication by $d$ identifies $\Z_k$ with $d\Z_n$.
Then, $g$ and all $f_r$ are permutations. Moreover,
$$
m\in D_n(f)\quad\Longleftrightarrow\quad
m\bmod d\in D_d(g)\ \text{and}\
m\bmod k\in\bigcap_{r=0}^{d-1}D_k(f_r).
$$
\end{lemma}

\begin{proof}
The map $g$ is well defined because $f$ preserves congruences modulo
$d$, and it is surjective because $f$ is. Hence, $g$ is a permutation.
Moreover, $f$ maps the residue class of $r$ modulo $d$ into the residue
class of $f(r)$. Both classes have $k$ elements, and $f$ is injective.
Consequently, this restriction is a bijection, and $f_r$ is a
well-defined permutation of $\Z_k$.

For $x=r+dt$, we have
$f(x)+mx=f(r)+mr+d(f_r(t)+mt)$. Thus, the map $x\mapsto f(x)+mx$
preserves congruences modulo $d$, and the map that it induces on $\Z_d$
is $y\mapsto g(y)+my$. On the class of $r$, its injectivity is equivalent
to that of $t\mapsto f_r(t)+mt$. If $f(x)+mx$ is a permutation, its
induced map on $\Z_d$ is surjective, and each restriction is injective.
Conversely, if the induced map is a permutation, distinct input classes
have distinct output classes. If every restriction is also injective,
then the whole map is injective. This proves the equivalence.
\end{proof}

We call $g$ the \emph{quotient permutation} of $f$, and we call the maps
$f_r$ the \emph{residue-class permutations} of $f$. By
Lemma~\ref{lem:split-general}, the good slopes of $f$ are determined by
those of these smaller permutations.

We close this section with two elementary facts, which are used in the
proof of Theorem~\ref{thm:universal}. The first one is a form of
Lagrange interpolation over commutative rings.

\begin{lemma}\label{lem:interp}
Let $A$ be a commutative ring with identity, and let $a_1,\dots,a_t\in A$
have pairwise invertible differences. If $P\in A[T]$ has degree less than
$t$ and $P(a_i)=0$ for $1\le i\le t$, then $P=0$.
\end{lemma}

\begin{proof}
Write $P(T)=\sum_{j=0}^{t-1}c_jT^j$. The conditions $P(a_i)=0$ form a
homogeneous linear system in $c_0,\dots,c_{t-1}$. Its matrix is the
Vandermonde matrix with rows $(1,a_i,\dots,a_i^{t-1})$, and its determinant
$\prod_{i<k}(a_k-a_i)$ is a unit of $A$. Hence, the matrix is invertible
over $A$, and all $c_j$ vanish.
\end{proof}

The hypothesis on the differences cannot be omitted. For instance, the
nonzero polynomial $3T\in\Z_9[T]$ vanishes at $0$, $3$ and $6$.

The second fact concerns power sums. Over finite fields, power sums are
a standard tool in the study of permutation polynomials, for instance
through Hermite's criterion \cite{LN}. We need the following congruences
over $\Z_n$.

\begin{lemma}\label{lem:powersum}
Let $p$ be a prime with $p^e\,\|\,n$, where $e\ge1$, and put
$s_i=\sum_{x=0}^{n-1}x^i$ for $i\ge1$.
\begin{enumerate}
\item If $p$ is odd, then
$s_i\equiv\frac np\sum_{x=0}^{p-1}x^i\pmod{p^e}$. Consequently,
$s_i\equiv0\pmod{p^e}$ for $1\le i\le p-2$, and
$s_{p-1}\equiv-n/p\pmod{p^e}$.
\item If $p=2$, then $s_1\equiv-n/2\pmod{2^e}$.
\end{enumerate}
In particular, $s_{p-1}\not\equiv0\pmod{p^e}$ in both cases.
\end{lemma}

\begin{proof}
Suppose first that $p$ is odd. For $k\ge1$, put
$S_k=\sum_{x=0}^{p^k-1}x^i$. We show by induction on $k$ that
$S_k\equiv p^{k-1}S_1\pmod{p^k}$. This is trivial for $k=1$. Let
$k\ge2$, and write $x=u+vp^{k-1}$ with $0\le u<p^{k-1}$ and $0\le v<p$.
Since $2(k-1)\ge k$, the binomial theorem gives
$(u+vp^{k-1})^i\equiv u^i+iu^{i-1}vp^{k-1}\pmod{p^k}$. We sum over $u$ and
$v$. The first terms contribute $pS_{k-1}$. As $p$ is odd, the sum
$\sum_{v=0}^{p-1}v=p(p-1)/2$ is divisible by $p$. Hence, the second terms
contribute a multiple of $p^k$, and $S_k\equiv pS_{k-1}\pmod{p^k}$. Now
the induction hypothesis $S_{k-1}\equiv p^{k-2}S_1\pmod{p^{k-1}}$ gives
$S_k\equiv p^{k-1}S_1\pmod{p^k}$.

Every residue modulo $p^e$ occurs exactly $n/p^e$ times among
$0,1,\dots,n-1$, and $x^i$ modulo $p^e$ depends only on $x$ modulo $p^e$.
Therefore,
$$
s_i\equiv\frac{n}{p^e}\,S_e\equiv\frac{n}{p^e}\,p^{e-1}S_1
=\frac np\,S_1\pmod{p^e}.
$$
Since $n/p$ is divisible by $p^{e-1}$, the residue of $(n/p)S_1$ modulo
$p^e$ depends only on $S_1$ modulo $p$. For $i=p-1$, Fermat's theorem gives
$S_1\equiv p-1\equiv-1\pmod p$. For $1\le i\le p-2$, the cyclic group
$\F_p^\times$ has order $p-1>i$. Hence, there is $t\in\F_p^\times$ with
$t^i\ne1$. Multiplication by $t$ permutes $\F_p$. Thus,
$t^iS_1\equiv S_1$, and hence $S_1\equiv0\pmod p$. This proves (i).

For $p=2$, we have $s_1=n(n-1)/2=-n/2+n^2/2$, and $n^2/2$ is divisible by
$2^{2e-1}$, hence by $2^e$.

Finally, $n/p$ is divisible by $p^{e-1}$ but not by $p^e$. Hence,
$-n/p\not\equiv0\pmod{p^e}$.
\end{proof}

\section{Universal bounds}\label{sec:universal}

In this section, we prove the bounds $|D_n(f)|\le\varphi(n)$ and
$\mu_n(f)\le p_0(n)-1$ stated in the introduction. Both follow from one
fact, namely that the set of good slopes misses a residue class modulo
every prime divisor of $n$. We prove this fact with the power sums of
Lemma~\ref{lem:powersum}, and we then derive two classical consequences.

\begin{theorem}\label{thm:universal}
Let $f$ be a permutation of $\Z_n$, and let $p$ be a prime divisor of $n$.
Then, $D_n(f)$ does not meet every residue class modulo $p$. Consequently,
the following hold.
\begin{enumerate}
\item There is $a\in\Z_n^\times$ such that $D_n(f)\subseteq\Z_n^\times-a$.
In particular, $|D_n(f)|\le\varphi(n)$.
\item $\mu_n(f)\le p_0(n)-1$.
\item If $|D_n(f)|=\varphi(n)$, then $D_n(f)=\Z_n^\times-a$ for the unit
$a$ in (i), and $\mu_n(f)=p_0(n)-1$.
\end{enumerate}
\end{theorem}

\begin{proof}
Let $p^e\,\|\,n$. We choose integer representatives for the values of $f$
and consider
$$
Q(T)=\sum_{x=0}^{n-1}\bigl(f(x)+xT\bigr)^{p-1}-\sum_{x=0}^{n-1}x^{p-1}
\in\Z[T].
$$
If $m\in D_n(f)$, then the integers $f(x)+mx$ with $0\le x<n$ represent
every residue class modulo $n$ exactly once. Hence, $Q(m)\equiv0\pmod n$,
and in particular $Q(m)\equiv0\pmod{p^e}$. The coefficient of $T^{p-1}$ in
$Q$ is $s_{p-1}=\sum_xx^{p-1}$, which is nonzero modulo $p^e$ by
Lemma~\ref{lem:powersum}. Suppose that $D_n(f)$ meets every residue class
modulo $p$, and choose $m_0,\dots,m_{p-1}\in D_n(f)$ with
$m_i\equiv i\pmod p$. Their pairwise differences are units of
$\Z_{p^e}$. The reduction of $Q$ modulo $p^e$ has degree at most $p-1$
and vanishes at these $p$ elements. By Lemma~\ref{lem:interp}, applied with
$A=\Z_{p^e}$, this reduction is zero. This contradicts
$s_{p-1}\not\equiv0\pmod{p^e}$ and proves the first assertion.

For every prime $p\mid n$, choose a residue class $c_p$ modulo $p$ that
contains no good slope. Since $0\in D_n(f)$, we have $c_p\neq0$. By the
Chinese remainder theorem, there is $a\in\Z_n$ with $a\equiv-c_p\pmod p$
for every prime $p\mid n$, and this $a$ is a unit. If $m\in D_n(f)$, then
$m\not\equiv c_p\equiv-a\pmod p$ for every $p$. Hence, $m+a$ is divisible
by no prime divisor of $n$. Thus, $m+a\in\Z_n^\times$, which proves
$D_n(f)\subseteq\Z_n^\times-a$. The bound in (i) follows from
$|\Z_n^\times-a|=\varphi(n)$.

The set $\Z_n^\times-a$ consists of the $m\in\Z_n$ with
$m\not\equiv-a\pmod p$ for every prime $p\mid n$. We apply
Lemma~\ref{lem:admissible}(ii) with $T_p$ equal to the set of residues
modulo $p$ different from $-a$. It shows that every subset of
$\Z_n^\times-a$ with pairwise invertible differences has at most
$\min_{p\mid n}(p-1)=p_0(n)-1$ elements. Every admissible subset of
$D_n(f)$ lies in $\Z_n^\times-a$, and this proves (ii). If
$|D_n(f)|=\varphi(n)$, then (i) gives $D_n(f)=\Z_n^\times-a$. In this case,
Lemma~\ref{lem:admissible}(ii) also provides an admissible subset with
$p_0(n)-1$ elements, which proves (iii).
\end{proof}

We have seen that an affine permutation $f(x)=ax+b$ satisfies
$D_n(f)=\Z_n^\times-a$. Thus, both bounds in Theorem~\ref{thm:universal}
are attained, and the question is how close non-affine permutations can
come to them. When $n=p$ is prime, the theorem only says that
$D_p(f)\neq\F_p$. For non-affine permutations of a prime field, the much
stronger bound of R\'edei is recalled in Theorem~\ref{thm:fields}
below.

\begin{corollary}\label{cor:classical}
Let $f$ be a permutation of $\Z_n$.
\begin{enumerate}
\item If $n$ is even, then every good slope of $f$ is even. In particular,
$x\mapsto f(x)+x$ is not a permutation.
\item If $x\mapsto f(x)+x$ and $x\mapsto f(x)-x$ are permutations, then
$\gcd(n,6)=1$.
\end{enumerate}
\end{corollary}

\begin{proof}
By Theorem~\ref{thm:universal}, $D_n(f)$ misses a residue class modulo
every prime divisor of $n$. In (i), this class is the odd one, since
$0\in D_n(f)$. In (ii), the set $\{-1,0,1\}\subseteq D_n(f)$ meets both
residue classes modulo $2$ and all three residue classes modulo $3$.
Hence, neither $2$ nor $3$ divides $n$.
\end{proof}

Recall that a complete mapping of $\Z_n$ is a permutation $f$ for which
$x\mapsto f(x)+x$ is also a permutation. Thus, part (i) is the classical
fact that a cyclic group of even order has no complete mapping, see
\cite{Evans}. Part (ii) is the necessary condition in
P\'olya's theorem on the toroidal $n$-queens problem
\cite{BellStevens,Polya}. Indeed, $n$ queens placed on the cells
$(x,f(x))$ of the $n\times n$ torus are pairwise non-attacking exactly when
$f$, $f(x)+x$ and $f(x)-x$ are permutations. Conversely, such a
permutation exists whenever $\gcd(n,6)=1$. In this case, $f(x)=2x$ is a
permutation, and so are $f(x)-x=x$ and $f(x)+x=3x$.

\section{Non-affine permutations attaining the bounds}
\label{sec:constructions}

In this section, we show that, for every composite $n\ge6$, non-affine
permutations attain both bounds of Theorem~\ref{thm:universal}. The
examples come from maps that are strictly triangular with respect to a
chain of ideals. For a prime modulus, the situation is different. There,
we use the direction theorems over finite fields, which give a much
smaller bound. Theorem~\ref{thm:maxima} combines both cases.

The ideals of $\Z_n$ are the sets $d\Z_n$ with $d\mid n$, and
$x-y\in d\Z_n$ means $x\equiv y\pmod d$.

\begin{definition}\label{def:triangular}
Let $\Z_n=I_0\supsetneq I_1\supsetneq\cdots\supsetneq I_r=\{0\}$ be a chain
of ideals. A map $g:\Z_n\to\Z_n$ is \emph{strictly triangular} with respect
to this chain if $x-y\in I_k$ implies $g(x)-g(y)\in I_{k+1}$, for all
$x,y\in\Z_n$ and $0\le k<r$. A permutation $f$ of $\Z_n$ is \emph{of
triangular type} if $f(x)=ax+b+g(x)$ for some $a\in\Z_n^\times$, some
$b\in\Z_n$ and some map $g$ that is strictly triangular with respect to
some chain of ideals.
\end{definition}

Constant maps are strictly triangular with respect to every chain. Hence,
every affine permutation is of triangular type.

\begin{proposition}\label{prop:triangular}
Let $g:\Z_n\to\Z_n$ be strictly triangular with respect to a chain of
ideals. Then, $x\mapsto g(x)+ux$ is a permutation of $\Z_n$ for every
$u\in\Z_n^\times$. Consequently, for all $a\in\Z_n^\times$ and
$b\in\Z_n$, the map $f(x)=ax+b+g(x)$ is a permutation with
$D_n(f)=\Z_n^\times-a$ and $\mu_n(f)=p_0(n)-1$.
\end{proposition}

\begin{proof}
Let the chain be $\Z_n=I_0\supsetneq\cdots\supsetneq I_r=\{0\}$, and let
$u\in\Z_n^\times$. Suppose that $g(x)+ux=g(y)+uy$ for some $x\neq y$, and
let $k$ be the largest index with $x-y\in I_k$. Then, $k<r$, and the
hypothesis gives $u(x-y)=g(y)-g(x)\in I_{k+1}$. Multiplying by $u^{-1}$, we
get $x-y\in I_{k+1}$, which contradicts the choice of $k$. Hence,
$x\mapsto g(x)+ux$ is injective, and therefore a permutation.

Now let $m\in\Z_n$ with $a+m\in\Z_n^\times$. Then,
$f(x)+mx=g(x)+(a+m)x+b$ is a permutation by the first part. Taking $m=0$,
we see that $f$ is a permutation. Thus, $\Z_n^\times-a\subseteq D_n(f)$,
and Theorem~\ref{thm:universal} gives $D_n(f)=\Z_n^\times-a$ and
$\mu_n(f)=p_0(n)-1$.
\end{proof}

\begin{example}\label{ex:primepower}
Let $p$ be a prime, let $e\ge2$, and put $n=p^e$. Consider the chain
$I_k=p^k\Z_{p^e}$, where
$0\le k\le e$. A map $g$ is strictly triangular with respect to it if
$x\equiv y\pmod{p^k}$ implies $g(x)\equiv g(y)\pmod{p^{k+1}}$ for
$0\le k<e$. Every map $g=pP$ with $P\in\Z[X]$ has this property, because
$x-y$ divides $P(x)-P(y)$. Taking $P=X^2$ and $a=1$ in
Proposition~\ref{prop:triangular}, we see that $f(x)=x+px^2$ is a
permutation of $\Z_{p^e}$ with
$$
D_{p^e}(f)=\{m\in\Z_{p^e}:m\not\equiv-1\pmod p\}
\quad\text{and}\quad \mu_{p^e}(f)=p-1.
$$
For instance, $\{0,1,\dots,p-2\}$ is an admissible set. Every affine map
$h$ satisfies $h(2)-2h(1)+h(0)=0$, whereas $f(2)-2f(1)+f(0)=2p$. This is
nonzero in $\Z_{p^e}$ unless $(p,e)=(2,2)$. Hence, $f$ is not affine when
$p$ is odd, and when $p=2$ and $e\ge3$.

For $n=9$, the permutation $f(x)=x+3x^2$ is given by
$$
\begin{array}{|c|ccccccccc|}
\hline
x&0&1&2&3&4&5&6&7&8\\ \hline
f(x)&0&4&5&3&7&8&6&1&2\\ \hline
\end{array}
$$
and $D_9(f)=\{0,1,3,4,6,7\}$. Thus, $f$ has six good slopes out of nine,
while $\mu_9(f)=2$. For $n=8$, the permutation $f(x)=x+2x^2$ of $\Z_8$
satisfies $\bigl(f(0),\dots,f(7)\bigr)=(0,3,2,5,4,7,6,1)$. Its good
slopes are $0$, $2$, $4$ and $6$, and $\mu_8(f)=1$.
\end{example}

\begin{figure}[htbp]
\centering
\begin{tikzpicture}[x=.49cm,y=.49cm,every node/.style={font=\small}]
\begin{scope}[shift={(0,0)}]
\foreach \i in {0,...,9}{
  \draw[black!22,line width=.25pt] ({\i-.5},-.5) -- ({\i-.5},8.5);
  \draw[black!22,line width=.25pt] (-.5,{\i-.5}) -- (8.5,{\i-.5});
}
\draw[black!60,line width=.5pt] (-.5,-.5) rectangle (8.5,8.5);
\foreach \i in {0,...,8}{
  \node[below] at (\i,-.65) {$\i$};
  \node[left] at (-.65,\i) {$\i$};
}
\node at (4,-2.05) {$x$};
\node[rotate=90] at (-2,4) {$f(x)+mx$};
\node at (4,9.65) {(a) $m=1$};
\foreach \x/\y in {0/0,1/5,2/7,3/6,4/2,5/4,6/3,7/8,8/1}{
  \fill (\x,\y) circle[radius=2.2pt];
}
\node at (4,-3.1) {$1\in D_9(f)$};
\end{scope}
\begin{scope}[shift={(14.5,0)}]
\foreach \row in {0,6}{
  \fill[black!10] (-.5,{\row-.5}) rectangle (8.5,{\row+.5});
}
\foreach \i in {0,...,9}{
  \draw[black!22,line width=.25pt] ({\i-.5},-.5) -- ({\i-.5},8.5);
  \draw[black!22,line width=.25pt] (-.5,{\i-.5}) -- (8.5,{\i-.5});
}
\draw[black!60,line width=.5pt] (-.5,-.5) rectangle (8.5,8.5);
\foreach \i in {0,...,8}{
  \node[below] at (\i,-.65) {$\i$};
  \node[left] at (-.65,\i) {$\i$};
}
\node at (4,-2.05) {$x$};
\node[rotate=90] at (-2,4) {$f(x)+mx$};
\node at (4,9.65) {(b) $m=2$};
\foreach \x/\y in {0/0,1/6,2/0,3/0,4/6,5/0,6/0,7/6,8/0}{
  \fill (\x,\y) circle[radius=2.2pt];
}
\node at (4,-3.1) {$2\notin D_9(f)$};
\end{scope}

\end{tikzpicture}
\caption{The graphs of $x\mapsto f(x)+mx$ for $f(x)=x+3x^2$ over
$\Z_9$. Coordinates are represented by $0,\dots,8$.
For $m=1$, all nine images are distinct. For $m=2$, the image is
$\{0,6\}$, and the two shaded rows contain six and three points,
respectively.}
\label{fig:good-slopes-nine}
\end{figure}
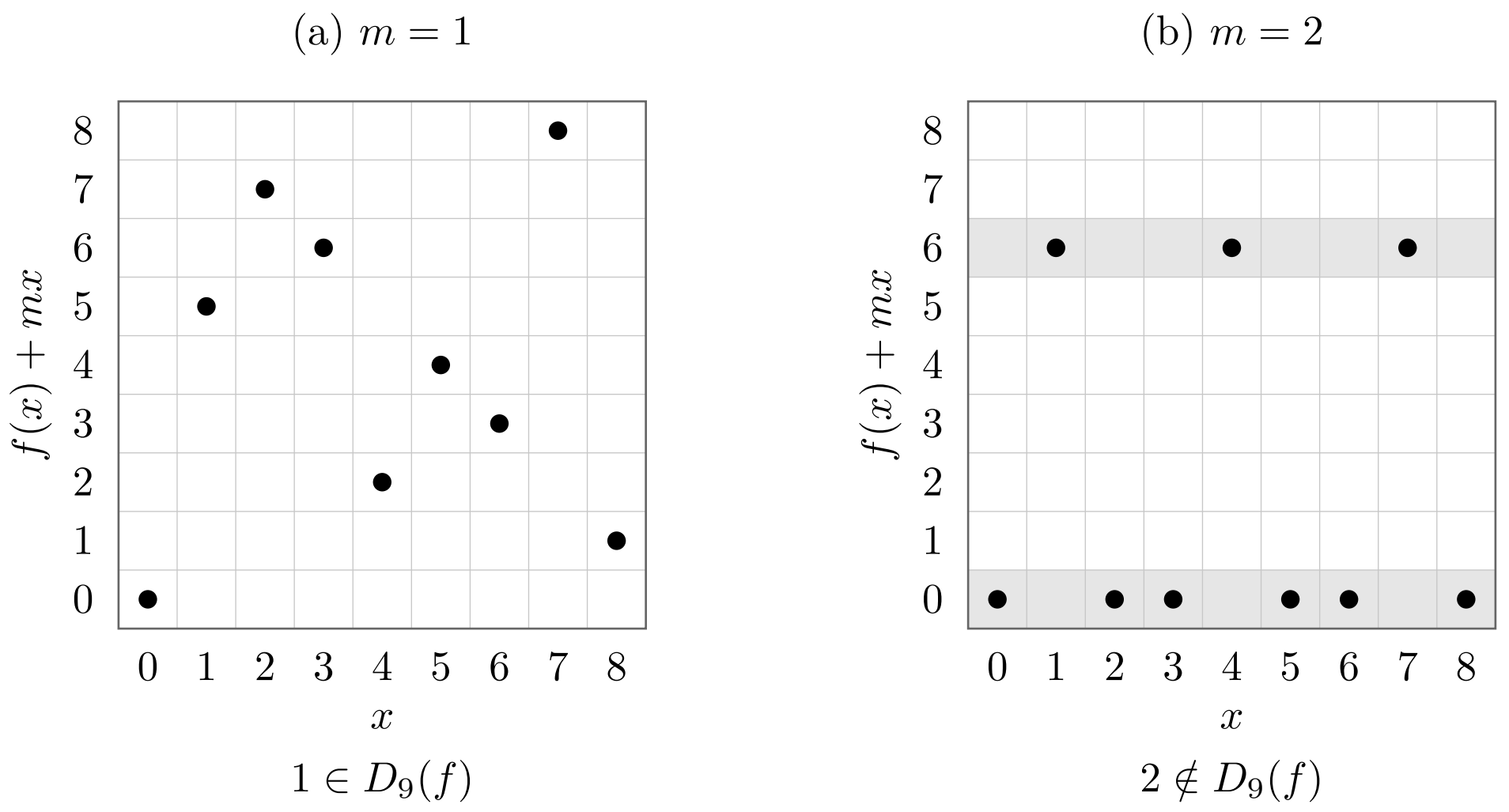

For the permutation $f(x)=x+3x^2$ of $\Z_9$ in
Example~\ref{ex:primepower}, Figure~\ref{fig:good-slopes-nine}
illustrates the definition of a good slope. Each column contains one
point. A slope is good precisely when each row also contains one point.

\begin{example}\label{ex:coprime}
Let $n=uv$ with $u,v\ge2$ and $\gcd(u,v)=1$. Let $\varepsilon\in\Z_n$ be
the idempotent with $\varepsilon\equiv1\pmod u$ and
$\varepsilon\equiv0\pmod v$. For a map $\theta:\Z_v\to\Z$, put
$g(x)=\varepsilon\,\theta(x\bmod v)$. Then, $g$ takes values in $v\Z_n$, and
it is constant on the cosets of $v\Z_n$. Hence, $g$ is strictly triangular
with respect to the chain $\Z_n\supsetneq v\Z_n\supsetneq\{0\}$. By
Proposition~\ref{prop:triangular}, $f(x)=x+g(x)$ is a permutation with
$D_n(f)=\Z_n^\times-1$ and $\mu_n(f)=p_0(n)-1$. In the coordinates
$(s,t)=(x\bmod u,\,x\bmod v)$, we have $f(s,t)=(s+\theta(t),t)$. If
$\theta$ is not constant modulo $u$, then $f$ is not affine. Indeed, the
first coordinate of an affine map $x\mapsto ax+b$ is $as+b$, which does
not depend on $t$. Every $n$ with at least two distinct prime divisors has
a factorisation of this kind.

For $n=6$, $u=2$, $v=3$ and $\theta(t)=[t=1]$, we have $\varepsilon=3$. The
map $f(x)=x+3\,[x\equiv1\pmod3]$ is the transposition $(1\ 4)$, with
$D_6(f)=\{0,4\}$ and $\mu_6(f)=1$. For $n=15$, $u=3$, $v=5$ and the same
$\theta$, we have $\varepsilon=10$. The map
$f(x)=x+10\,[x\equiv1\pmod5]$ is the $3$-cycle $(1\ 11\ 6)$, and
$D_{15}(f)=\Z_{15}^\times-1=\{0,1,3,6,7,10,12,13\}$. Thus,
$|D_{15}(f)|=8>15/2$ and $\mu_{15}(f)=2$.
\end{example}

For prime $n$, we need the theory of directions over finite fields. Let
$q$ be a prime power. The lines of $\F_q^2$ have $q+1$ directions,
namely the slopes $a\in\F_q$ of
the lines $\{(x,y):y=ax+c\}$ and the vertical direction. A set of points
\emph{determines} a direction if some line with this direction contains at
least two of its points. For a map $h:\F_q\to\F_q$, which need not be a
permutation, we write $D_{\F_q}(h)$ for the set of $m\in\F_q$ such that
$x\mapsto h(x)+mx$ is a permutation of $\F_q$. When $q=p$ is prime and $h$
is a permutation, this is the set $D_p(h)$ of Definition~\ref{def:slopes}.
The graph of $h$ never determines the vertical direction. As explained in
the introduction, it determines the slope $-m$ exactly when
$m\notin D_{\F_q}(h)$. We use the following theorem of Ball, which
completed the work of \cite{BBBSS}.

\begin{theorem}[{\cite[Theorem~1.1]{Ball}}]\label{thm:ball}
Let $q=p^t$ with $p$ prime and $t\ge1$, let $h:\F_q\to\F_q$ be a map,
and let $N$ be the number of directions determined by the graph of $h$. Let $s$ be the
largest power of $p$ such that every line with a determined direction
meets the graph of $h$ in $0$ or in a multiple of $s$ points. Then, one of
the following holds.
\begin{enumerate}
\item $s=1$ and $N\ge(q+3)/2$.
\item $1<s<q$, $\F_s$ is a subfield of $\F_q$, and
$q/s+1\le N\le(q-1)/(s-1)$.
\item $s=q$ and $N=1$.
\end{enumerate}
Moreover, if $s>2$, then the map $x\mapsto h(x)-h(0)$ is
$\F_s$-linear.
\end{theorem}

Over a finite field, Theorem~\ref{thm:ball} gives the following bounds
for good slopes.

\begin{theorem}\label{thm:fields}
Let $q$ be a prime power and let $h:\F_q\to\F_q$ be a map.
\begin{enumerate}
\item If $q$ is odd and $h$ is not the sum of an additive map and a
constant, then $|D_{\F_q}(h)|\le(q-3)/2$.
\item If $|D_{\F_q}(h)|>q/2$, then $h$ is the sum of an additive map and a
constant.
\end{enumerate}
\end{theorem}

\begin{proof}
Write $q=p^t$, where $p$ is prime and $t\ge1$. For $q$ prime, the
additive maps of $\F_q$ are linear, and (i) is R\'edei's theorem
\cite{Redei}. In general, we deduce both parts from
Theorem~\ref{thm:ball}. Since the vertical direction is not determined,
we have $N=q-|D_{\F_q}(h)|$. If $s>2$, then $h$ is the sum of an additive
map and the constant $h(0)$.

Suppose that $q$ is odd and that $h$ is not the sum of an additive map and
a constant. Then, $s\le2$. As $s$ is a power of the odd prime $p$, we get
$s=1$. Hence, $N\ge(q+3)/2$, and (i) follows.

Now suppose that $|D_{\F_q}(h)|>q/2$. Then, $q>2$, because a map of $\F_2$
has at most one good slope. Moreover, $N<q/2$. This excludes $s=1$. It
also excludes $s=2$, since then $s<q$ and $N\ge q/2+1$. Therefore, $s>2$,
and (ii) follows.
\end{proof}

We also need a form of R\'edei's theorem for sets of points that are not
given as graphs. Let $p$ be a prime. From now on, we write the directions
of $\F_p^2$ as one-dimensional subspaces $\langle v\rangle$. The lines with
direction $\langle v\rangle$ are the cosets $P+\langle v\rangle$, and the
direction $\langle(0,1)\rangle$ is \emph{vertical}. Every other direction is
$\langle(1,-m)\rangle$ for a unique $m\in\F_p$, and two points $(r,u)$ and
$(r',u')$ lie on a common line with this direction if and only if
$u+mr=u'+mr'$. A set $S\subseteq\F_p^2$ of $p$ points is \emph{injective
in the direction} $\langle v\rangle$ if every line with direction
$\langle v\rangle$ contains exactly one point of $S$. A set of $p$ points
that is contained in a line is equal to this line. Moreover, a set of
$p$ points is injective in a direction exactly when it does not
determine this direction. Thus, part (ii) of the next lemma is the
theorem of R\'edei and Megyesi, which states that a set of $p$ points
of $\F_p^2$ that is not a line determines at least $(p+3)/2$ directions
\cite{Redei}. Somlai has recently given a new proof of this theorem
\cite{Somlai}.

\begin{lemma}\label{lem:pointsets}
Let $p$ be an odd prime, and let $S\subseteq\F_p^2$ be a set of $p$ points
that is not a line.
\begin{enumerate}
\item If $S$ meets every vertical line in exactly one point, then $S$ is
injective in at most $(p-3)/2$ non-vertical directions.
\item The set $S$ is injective in at most $(p-1)/2$ directions.
\end{enumerate}
\end{lemma}

\begin{proof}
In (i), $S$ is the graph $\{(r,h(r)):r\in\F_p\}$ of a map
$h:\F_p\to\F_p$. By the description of the lines above, $S$ is injective
in the direction $\langle(1,-m)\rangle$ if and only if
$r\mapsto h(r)+mr$ is a permutation, that is, if and only if
$m\in D_{\F_p}(h)$. The graph of an affine map is a line. Hence, $h$ is
not affine. Theorem~\ref{thm:fields}(i) gives (i).

For (ii), we may assume that $S$ is injective in some direction
$\langle v\rangle$. Choose $A\in\mathrm{GL}_2(\F_p)$ with
$Av\in\langle(0,1)\rangle$. The map $A$ permutes the lines and the
directions of $\F_p^2$, and $AS$ is not a line. Moreover, $AS$ is
injective in the vertical direction, which means that it meets every
vertical line in exactly one point. By (i), $AS$ is injective in at most
$(p-3)/2$ non-vertical directions. Together with the vertical direction,
this gives at most $(p-1)/2$ directions for $AS$, and hence for $S$.
\end{proof}

The next proposition treats the index-two case of a construction attributed
to Megyesi in \cite{BBS}, see also \cite[Example~15]{BS}. We include the short
proof because we need the explicit set of good slopes.

\begin{proposition}\label{prop:megyesi}
Let $q\ge5$ be an odd prime power, let $H$ be the group of nonzero squares
of $\F_q$, and let
$$
h(x)=\begin{cases}x&\text{if }x\in H,\\0&\text{if }x\notin H.\end{cases}
$$
Then, $h$ is not the sum of an additive map and a constant, and
$$
D_{\F_q}(h)=\bigl\{(z-1)^{-1}:z\in H\setminus\{1\}\bigr\},\quad
|D_{\F_q}(h)|=\frac{q-3}{2}.
$$
For every $c\in D_{\F_q}(h)$, the map $f(x)=h(x)+cx$ is a permutation of
$\F_q$ with $D_{\F_q}(f)=D_{\F_q}(h)-c$, and $f$ is not the sum of an
additive map and a constant.
\end{proposition}

\begin{proof}
Let $b\in\F_q$. The map $x\mapsto h(x)+bx$ fixes $0$, multiplies the
elements of $H$ by $b+1$, and multiplies the non-squares by $b$. If $b=0$,
it maps every non-square to $0$. If $b=-1$, it maps $H$ to $0$. In these
two cases, it is not a permutation. Otherwise, its nonzero values form the
two cosets $(b+1)H$ and $bH'$ of $H$, where $H'=\F_q^\times\setminus H$.
Hence, it is a permutation exactly when these two cosets are different. As
$H$ has index $2$ in $\F_q^\times$ and $bH'\neq bH$, this happens exactly
when $(b+1)H=bH$, that is, when $z=1+b^{-1}\in H$. The conditions $b\ne0$
and $b\neq-1$ mean $z\neq1$ and $z\neq0$, and then $b=(z-1)^{-1}$. This
gives the stated description of $D_{\F_q}(h)$, whose size is
$|H|-1=(q-3)/2$.

Suppose that $h$ is the sum of an additive map and a constant. Since
$h(0)=0$, the constant is $0$, and $h$ is additive. Let $z$ be a
non-square. Then, $h(z)=0$ and $h(1)=1$. By additivity, $h(z+1)=1$. But
$h(z+1)\in\{0,z+1\}$ and $z+1\neq1$, a contradiction. The map $f$ differs
from $h$ by the linear map $x\mapsto cx$. Hence, $f$ is not the sum of an
additive map and a constant either. Moreover, $f(x)+mx=h(x)+(c+m)x$.
Hence, $D_{\F_q}(f)=D_{\F_q}(h)-c$. This set contains $0$, and therefore $f$ is a
permutation. Finally, $D_{\F_q}(h)$ is nonempty because $q\ge5$.
\end{proof}

\begin{example}\label{ex:p7}
For $q=7$, we have $H=\{1,2,4\}$ and
$\bigl(h(0),\dots,h(6)\bigr)=(0,1,2,0,4,0,0)$.
Proposition~\ref{prop:megyesi} gives
$D_{\F_7}(h)=\{1^{-1},3^{-1}\}=\{1,5\}$. With $c=1$, the permutation
$f(x)=h(x)+x$ is the $3$-cycle $(1\ 2\ 4)$, and $D_7(f)=\{0,4\}$ has
$(7-3)/2=2$ elements.
\end{example}

Proposition~\ref{prop:megyesi} shows that the bound of
Theorem~\ref{thm:fields}(i) is sharp for every odd $q\ge5$. For a prime
$q$, Lov\'asz and Schrijver proved that the graph of every
map $\F_q\to\F_q$ with exactly $(q-3)/2$ good slopes is the image of
the graph of $x\mapsto x^{(q+1)/2}$ under an affine transformation of
$\F_q^2$ \cite{LS}. The map $h$ of Proposition~\ref{prop:megyesi}
satisfies $h(x)=\bigl(x+x^{(q+1)/2}\bigr)/2$ for every $x\in\F_q$.
Hence, its graph is the image of the graph of $x\mapsto x^{(q+1)/2}$
under the linear map $(x,y)\mapsto\bigl(x,(x+y)/2\bigr)$.

We can now determine the largest values of $|D_n(f)|$ and $\mu_n(f)$ for
non-affine permutations. For $n\ge4$, let $A(n)$ and $B(n)$ denote these
largest values over all non-affine permutations $f$ of $\Z_n$. For $n=2$
and $n=3$, there are no non-affine permutations, since the number
$n\varphi(n)$ of affine permutations equals $n!$.

\begin{theorem}\label{thm:maxima}
We have $A(4)=B(4)=1$, and $A(p)=B(p)=(p-3)/2$ for every prime $p\ge5$.
For every composite $n\ge6$, we have $A(n)=\varphi(n)$ and
$B(n)=p_0(n)-1$, and some non-affine permutation $f$ with
$D_n(f)=\Z_n^\times-1$ attains both values.
\end{theorem}

\begin{proof}
Let $n\ge6$ be composite. By Theorem~\ref{thm:universal},
$A(n)\le\varphi(n)$ and $B(n)\le p_0(n)-1$. If $n$ is a prime power $p^e$,
then $e\ge2$ and $(p,e)\neq(2,2)$. In this case, Example~\ref{ex:primepower}
gives a non-affine permutation $f$ with $D_n(f)=\Z_n^\times-1$ and
$\mu_n(f)=p_0(n)-1$. Otherwise, $n$ has two distinct prime divisors, and
Example~\ref{ex:coprime} with $\theta(t)=[t=1]$ gives such a permutation.

Let $n=4$. By Lemma~\ref{lem:normal}, we may assume that $f(0)=0$. By
Corollary~\ref{cor:classical}(i), every good slope is even. Thus,
$D_4(f)\subseteq\{0,2\}$. The permutations of $\Z_4$ fixing $0$ are the
affine maps $x$ and $3x$, together with the four non-affine maps
$(0,1,3,2)$, $(0,2,1,3)$, $(0,2,3,1)$ and $(0,3,1,2)$, written as lists of
values. Adding $2x$, that is $(0,2,0,2)$, to them gives $(0,3,3,0)$,
$(0,0,1,1)$, $(0,0,3,3)$ and $(0,1,1,0)$, and none of these is a
permutation. Hence, every non-affine permutation of $\Z_4$ has exactly one
good slope, and $A(4)=B(4)=1$.

Let $p\ge5$ be prime. Theorem~\ref{thm:fields}(i) gives $A(p)\le(p-3)/2$,
and Proposition~\ref{prop:megyesi} gives equality. Finally, $B(p)=A(p)$,
because $\mu_p(f)=|D_p(f)|$ for every permutation $f$ of $\Z_p$.
\end{proof}

\begin{corollary}\label{cor:half}
There is a non-affine permutation $f$ of $\Z_n$ with $|D_n(f)|>n/2$ if and
only if $n$ is composite and $\varphi(n)>n/2$. In particular, this holds
for every odd prime power $p^e$ with $e\ge2$, and it fails for every even
$n$.
\end{corollary}

\begin{proof}
For $n\in\{2,3\}$, there are no non-affine permutations. For $n=4$ and for
primes $p\ge5$, Theorem~\ref{thm:maxima} gives $A(n)<n/2$. For composite
$n\ge6$, it gives $A(n)=\varphi(n)$. Finally,
$\varphi(n)=n\prod_{p\mid n}(1-1/p)\le n/2$ when $n$ is even, and
$\varphi(p^e)=p^e(1-1/p)>p^e/2$ when $p$ is odd.
\end{proof}

Corollary~\ref{cor:half} should be compared with
Theorem~\ref{thm:fields}(ii). Over a finite field, a map for which more
than half of all slopes are good is additive up to a constant. Over $\Z_n$
with $n$ composite, the same threshold does not force affinity.
Oddness of $n$ alone does not suffice for the existence statement, since
$\varphi(105)=48<105/2$.

\section{Extremal permutations of prime-power order}\label{sec:primepowers}

In this section, we classify the extremal permutations of prime-power
order. We show that they are exactly the permutations of triangular type,
describe them through their base-$p$ digits, and count them. The main
step is Lemma~\ref{lem:reduction}, which passes from $\Z_{p^e}$ to
$\Z_{p^{e-1}}$ inside every residue class modulo $p$.

Throughout this section, $p$ is a prime, $e\ge1$, and $n=p^e$. For a
nonzero $z\in\Z_{p^e}$, let
$\nu(z)$ be the largest integer $v<e$ such that $p^v$ divides $z$. We write
the elements of $\Z_{p^e}$ in base $p$ as
$x=x_0+x_1p+\dots+x_{e-1}p^{e-1}$ with digits $x_j\in\{0,\dots,p-1\}$, and
we identify the digits with elements of $\F_p$. For $x\neq y$, the integer
$\nu(x-y)$ is the first position in which the digits of $x$ and $y$
differ. By Theorem~\ref{thm:universal}(iii), a permutation $f$ of
$\Z_{p^e}$ with $\varphi(p^e)$ good slopes satisfies
$D_{p^e}(f)=\Z_{p^e}^\times-a$ for some unit $a$. This set consists of the
$m$ with $m\not\equiv-a\pmod p$. Hence, it determines $a$ modulo $p$. For
$e\ge2$, the map $t\mapsto pt$ is a bijection from $\Z_{p^{e-1}}$ onto
$p\Z_{p^e}$, and we use it to identify these two groups.

\begin{lemma}\label{lem:reduction}
Let $f$ be a permutation of $\Z_{p^e}$ with $D_{p^e}(f)=\Z_{p^e}^\times-a$,
where $a$ is a unit.
\begin{enumerate}
\item For all $x,y\in\Z_{p^e}$, we have $f(x)-f(y)\equiv a(x-y)\pmod p$.
\item Let $e\ge2$ and $r\in\Z_{p^e}$. There is a unique map
$f_r:\Z_{p^{e-1}}\to\Z_{p^{e-1}}$ with $f(r+pt)-f(r)=p\,f_r(t)$ for all
$t\in\Z_{p^{e-1}}$. Moreover, $f_r$ is a permutation with
$D_{p^{e-1}}(f_r)=\Z_{p^{e-1}}^\times-a$, where $a$ also denotes the
residue of $a$ modulo $p^{e-1}$.
\end{enumerate}
\end{lemma}

\begin{proof}
Let $x,y\in\Z_{p^e}$ with $x\not\equiv y\pmod p$. Then, $x-y$ is a unit, and
$m=-(f(x)-f(y))(x-y)^{-1}$ is the unique element with $f(x)+mx=f(y)+my$.
This $m$ is not a good slope. Hence, $m+a$ is not a unit, that is,
$m\equiv-a\pmod p$. This gives
$f(x)-f(y)=-m(x-y)\equiv a(x-y)\pmod p$. If $x\equiv y\pmod p$, we choose
$z\not\equiv x\pmod p$ and subtract the congruences for the pairs $(x,z)$
and $(y,z)$. This proves (i).

For (ii), part (i) shows that $f(r+pt)-f(r)\in p\Z_{p^e}$ for every $t$.
This gives the existence and the uniqueness of $f_r$. It is injective
because $f$ is, and hence it is a permutation. Now let $m\in\Z_{p^e}$ with
$m+a\in\Z_{p^e}^\times$, and put $g(x)=f(x)+mx$. Then, $g$ is a permutation,
and $g(r+pt)-g(r)=p\bigl(f_r(t)+mt\bigr)$, where $mt$ is computed in
$\Z_{p^{e-1}}$. As $g$ is injective, so is $t\mapsto f_r(t)+mt$. The
residues modulo $p^{e-1}$ of these elements $m$ are exactly the elements
of $\Z_{p^{e-1}}^\times-a$. Therefore,
$\Z_{p^{e-1}}^\times-a\subseteq D_{p^{e-1}}(f_r)$, and
Theorem~\ref{thm:universal} gives equality.
\end{proof}

We can now describe the extremal permutations of $\Z_{p^e}$.

\begin{theorem}\label{thm:classification}
Let $p$ be a prime, let $e\ge1$, and let $f$ be a permutation of
$\Z_{p^e}$. Write
$f(x)=y_0+y_1p+\dots+y_{e-1}p^{e-1}$ in base $p$. The following are
equivalent.
\begin{enumerate}
\item $|D_{p^e}(f)|=\varphi(p^e)$.
\item There is $a\in\Z_{p^e}^\times$ such that
$f(x)-f(y)\equiv a(x-y)\pmod{p^{\nu(x-y)+1}}$ for all distinct $x,y$.
\item There are $\alpha\in\F_p^\times$, $\beta\in\F_p$ and maps
$h_j:\F_p^{\,j}\to\F_p$ for $1\le j<e$ such that, in $\F_p$,
$$
y_0=\alpha x_0+\beta\quad\text{and}\quad
y_j=\alpha x_j+h_j(x_0,\dots,x_{j-1})
$$
for $1\le j<e$.
\item $f$ is of triangular type.
\end{enumerate}
In this case, $D_{p^e}(f)=\{m\in\Z_{p^e}:m\not\equiv-a\pmod p\}$ and
$\alpha\equiv a\pmod p$.
\end{theorem}

\begin{proof}
Suppose first that (i) holds. By Theorem~\ref{thm:universal}(iii), there
is a unit $a$ with $D_{p^e}(f)=\Z_{p^e}^\times-a$. We show by induction on
$e$ that (ii) holds for this $a$. If $\nu(x-y)=0$, then (ii) is the
congruence of Lemma~\ref{lem:reduction}(i). This settles $e=1$. Let
$e\ge2$, and let $x\neq y$ with $\nu(x-y)\ge1$. Then, $x\equiv y\pmod p$, and
we can write $x=r+pt$ and $y=r+pt'$ with $r\in\{0,\dots,p-1\}$ and distinct
$t,t'\in\Z_{p^{e-1}}$. Let $\nu'$ be defined on $\Z_{p^{e-1}}$ as $\nu$ is
defined on $\Z_{p^e}$. By Lemma~\ref{lem:reduction}(ii) and the induction
hypothesis, applied to $f_r$, we have
$f_r(t)-f_r(t')\equiv a(t-t')\pmod{p^{\nu'(t-t')+1}}$. Moreover, $x-y=p(t-t')$,
$f(x)-f(y)=p\bigl(f_r(t)-f_r(t')\bigr)$ and $\nu(x-y)=\nu'(t-t')+1$.
Multiplying the last congruence by $p$, we obtain (ii) for $x$ and $y$.

Next, (ii) implies (iv). Indeed, put $g(x)=f(x)-ax$. If $x\ne y$ and
$x-y\in p^k\Z_{p^e}$, then $\nu(x-y)\ge k$, and (ii) gives
$g(x)-g(y)\in p^{\nu(x-y)+1}\Z_{p^e}\subseteq p^{k+1}\Z_{p^e}$. Thus, $g$ is
strictly triangular with respect to the chain of ideals $p^k\Z_{p^e}$ with
$0\le k\le e$, and $f(x)=ax+g(x)$ is of triangular type. By
Proposition~\ref{prop:triangular}, (iv) implies (i), and
$D_{p^e}(f)=\Z_{p^e}^\times-a$ for the unit $a$ in (ii).

It remains to show that (ii) and (iii) are equivalent. Assume (ii). If
$x\equiv y\pmod{p^k}$ and $x\neq y$, then $\nu(x-y)\ge k$, and (ii) gives
$f(x)-f(y)\equiv a(x-y)\equiv0\pmod{p^k}$. Thus, for $0\le j<e$, the
residue of $f(x)$ modulo $p^{j+1}$ depends only on $x$ modulo $p^{j+1}$.
In other words, the digit $y_j$ depends only on $x_0,\dots,x_j$. Now fix
$j$ and $x_0,\dots,x_{j-1}$. Let $x$ and $x'$ have
these lower digits, different $j$th digits and equal higher digits, and
let $y'_j$ be the $j$th digit of $f(x')$. Then, $x-x'=p^j(x_j-x'_j)$ and
$\nu(x-x')=j$. Moreover, $f(x)$ and $f(x')$ have the same lower $j$ digits,
and hence $f(x)-f(x')\equiv p^j(y_j-y'_j)\pmod{p^{j+1}}$. Comparing with (ii), we
obtain $y_j-y'_j=\alpha(x_j-x'_j)$ in $\F_p$, where $\alpha\equiv a\pmod p$.
Hence, $y_j-\alpha x_j$ depends only on $x_0,\dots,x_{j-1}$. For $j=0$,
this means that $\beta=y_0-\alpha x_0$ is a constant. This proves (iii).

Conversely, assume (iii), and let $a\in\Z_{p^e}^\times$ with
$a\equiv\alpha\pmod p$. Let $x\neq x'$, write their digits as $x_j$ and
$x'_j$, and put $v=\nu(x-x')$. Then,
$x-x'\equiv p^v(x_v-x'_v)\pmod{p^{v+1}}$. By (iii), the digits of $f(x)$
and $f(x')$ agree below position $v$. In position $v$, they differ by
$\alpha(x_v-x'_v)$. Hence,
$f(x)-f(x')\equiv p^v\alpha(x_v-x'_v)\equiv a(x-x')\pmod{p^{v+1}}$,
which is (ii).
\end{proof}

Every choice of $\alpha\in\F_p^\times$, $\beta\in\F_p$ and maps
$h_1,\dots,h_{e-1}$ defines a permutation of $\Z_{p^e}$ through the
formulas in Theorem~\ref{thm:classification}(iii). Indeed, the input digits
can be recovered one at a time from the output digits.

\begin{corollary}\label{cor:count}
Let $p$ be a prime and $e\ge1$, and put $J=1+p+\dots+p^{e-1}$. The
number of permutations of $\Z_{p^e}$ with
$\varphi(p^e)$ good slopes is $(p-1)p^J$, and exactly $(p-1)p^{2e-1}$ of
them are affine. In particular, all of them are affine if and only if
$e=1$ or $(p,e)=(2,2)$.
\end{corollary}

\begin{proof}
The data in Theorem~\ref{thm:classification}(iii) are uniquely determined
by $f$. Indeed, $\alpha$ is determined by $D_{p^e}(f)$, and then $\beta$
and the maps $h_j$ are determined by the digits of $f$. By the remark
above, every choice of data gives a permutation with $\varphi(p^e)$ good
slopes. There are $p-1$ choices for $\alpha$, $p$ choices for $\beta$ and
$p^{p^j}$ choices for $h_j$. Hence, the number of such permutations is
$(p-1)\,p\prod_{j=1}^{e-1}p^{p^j}=(p-1)p^J$. Every affine permutation has
$\varphi(p^e)$ good slopes, and there are
$p^e\varphi(p^e)=(p-1)p^{2e-1}$ of them. The two numbers agree if and only
if $J=2e-1$. This holds for $e=1$. For $e=2$, it means $p+1=3$, that is,
$p=2$. For $e\ge3$, we have $J\ge1+2+\dots+2^{e-1}=2^e-1>2e-1$.
\end{proof}

\begin{corollary}\label{cor:param}
Let $p$ be a prime and $e\ge2$. A permutation $f$ of $\Z_{p^e}$ has
$\varphi(p^e)$ good slopes
if and only if $f(x)=ax+b+p\,H(x\bmod p^{e-1})$ for some
$a\in\{1,\dots,p-1\}$, some $b\in\Z_{p^e}$ and some
congruence-preserving map $H:\Z_{p^{e-1}}\to\Z_{p^{e-1}}$ with $H(0)=0$.
These data are uniquely determined by $f$. For $e=2$, the map
$H:\F_p\to\F_p$ is arbitrary, subject to $H(0)=0$.
\end{corollary}

\begin{proof}
Let $f$ have $\varphi(p^e)$ good slopes. Choose $a\in\{1,\dots,p-1\}$ with
$D_{p^e}(f)=\Z_{p^e}^\times-a$, put $b=f(0)$ and $G(x)=f(x)-ax-b$. By
Theorem~\ref{thm:classification}(ii), we have
$G(x)\equiv G(y)\pmod{p^{\nu(x-y)+1}}$ for $x\neq y$. Taking $y=0$, we see
that $G$ takes values in $p\Z_{p^e}$, because $G(0)=0$. If $x\neq y$ and
$x\equiv y\pmod{p^{e-1}}$, then $\nu(x-y)=e-1$, and hence $G(x)=G(y)$.
Thus, $G(x)=pH(x\bmod p^{e-1})$ for a unique map $H$ with $H(0)=0$. Let
$0\le k\le e-1$, let $u,v\in\Z_{p^{e-1}}$ be distinct with
$u\equiv v\pmod{p^k}$, and let $x,y\in\Z_{p^e}$ reduce to $u,v$. Then,
$\nu(x-y)\ge k$. Hence,
$p\bigl(H(u)-H(v)\bigr)=G(x)-G(y)\in p^{k+1}\Z_{p^e}$.
That is, $H(u)\equiv H(v)\pmod{p^k}$, and $H$ is congruence-preserving.

Conversely, let $H$ be congruence-preserving with $H(0)=0$, and put
$G(x)=pH(x\bmod p^{e-1})$. If $x\equiv y\pmod{p^k}$ with $0\le k<e$, then
the residues of $x$ and $y$ modulo $p^{e-1}$ are congruent modulo $p^k$.
Hence, $G(x)\equiv G(y)\pmod{p^{k+1}}$. Thus, $G$ is strictly triangular
with respect to the chain $p^k\Z_{p^e}$, and
Proposition~\ref{prop:triangular} shows that $f$ is a permutation with
$\varphi(p^e)$ good slopes. The data are unique because $a$ is determined
by $D_{p^e}(f)$, then $b=f(0)$, and then $H$ is determined by $f-ax-b$.
For $e=2$, the domain of $H$ is $\Z_p$, and the congruence condition
always holds.
\end{proof}

Corollary~\ref{cor:param} has a $p$-adic counterpart. A map $v$ of the
ring of $p$-adic integers is \emph{compatible} if $x\equiv y\pmod{p^k}$
implies $v(x)\equiv v(y)\pmod{p^k}$ for every $k\ge1$. Anashin proved
that, for every compatible map $v$ and all $p$-adic integers $a$ and
$b$ with $a\not\equiv0\pmod p$, the map $x\mapsto b+ax+pv(x)$ is
bijective modulo every power of $p$ \cite{Anashin}. For $e\ge2$,
Corollary~\ref{cor:param} shows that the extremal permutations of
$\Z_{p^e}$ are exactly the reductions modulo $p^e$ of these maps.

\begin{example}\label{ex:Z9}
For $n=9$, Corollary~\ref{cor:param} describes the permutations with six
good slopes as $f(x)=ax+b+3H(x\bmod3)$ with $a\in\{1,2\}$, $b\in\Z_9$ and
$H:\F_3\to\F_3$ with $H(0)=0$. There are $2\cdot9\cdot9=162=(3-1)3^{1+3}$
of them, in accordance with Corollary~\ref{cor:count}. Such an $f$ is
affine exactly when $H$ is linear. Hence, $108$ of them are not affine, and
$12$ of these fix $0$. For example, $a=1$, $b=0$ and $H(x)=[x=2]$ give the
permutation
$$
\begin{array}{|c|ccccccccc|}
\hline
x&0&1&2&3&4&5&6&7&8\\ \hline
f(x)&0&1&5&3&4&8&6&7&2\\ \hline
\end{array}
$$
with $D_9(f)=\{0,1,3,4,6,7\}$. In base $3$, this map is given by $y_0=x_0$
and $y_1=x_1+[x_0=2]$, as in Theorem~\ref{thm:classification}(iii). For
$n=27$, there are $2\cdot3^{13}=3188646$ permutations with $18$ good
slopes, and only $486$ of them are affine.
\end{example}

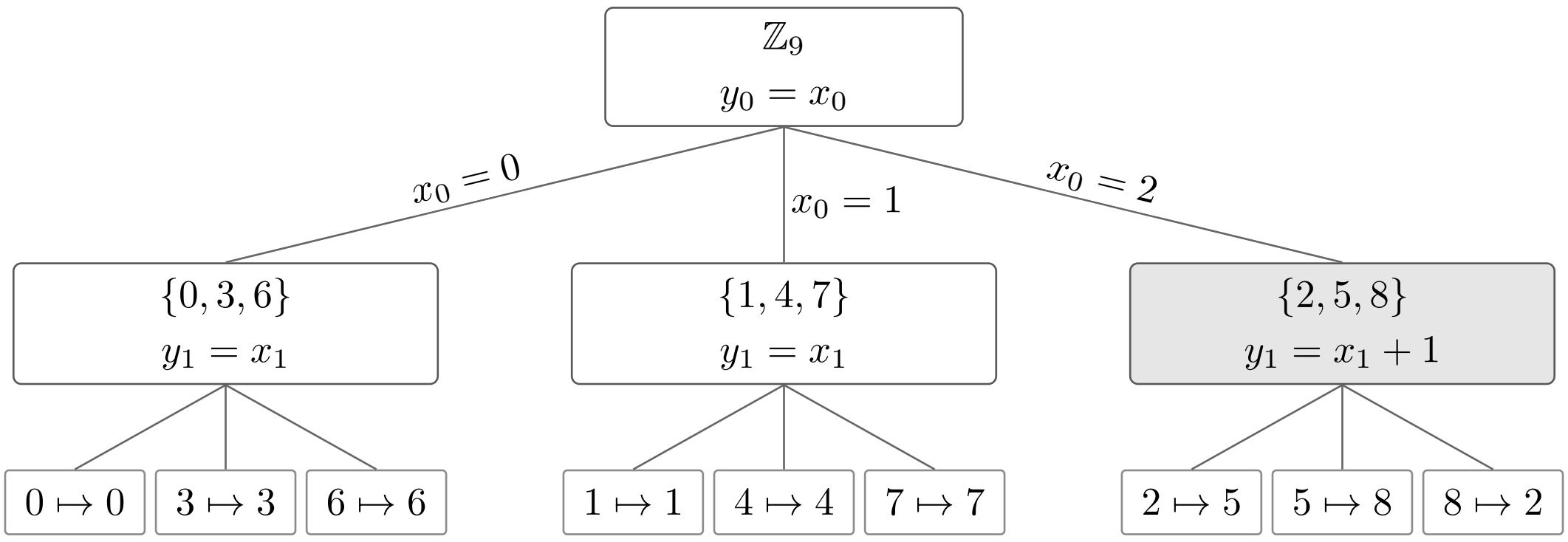
\begin{figure}[htbp]
\centering
\begin{tikzpicture}[
  every node/.style={font=\small},
  branch/.style={draw=black!65,rounded corners=2pt,align=center,
    minimum width=3.8cm,minimum height=1.05cm},
  leaf/.style={draw=black!45,rounded corners=1pt,fill=white,
    inner xsep=5pt,inner ysep=5pt},
  every path/.style={draw=black!60,line width=.5pt}]
\node[branch,minimum width=3.2cm,minimum height=.9cm] (root) at (0,3.9)
  {$\Z_9$\\[2pt]$y_0=x_0$};
\node[branch,fill=white] (b0) at (-5,1.6)
  {$\{0,3,6\}$\\[2pt]$y_1=x_1$};
\draw (root.south) -- node[above,sloped,pos=.56,fill=white,inner sep=1.5pt]
  {$x_0=0$} (b0.north);
\node[leaf] (leaf00) at (-6.35,0) {$0\mapsto 0$};
\draw (b0.south) -- (leaf00.north);
\node[leaf] (leaf01) at (-5,0) {$3\mapsto 3$};
\draw (b0.south) -- (leaf01.north);
\node[leaf] (leaf02) at (-3.65,0) {$6\mapsto 6$};
\draw (b0.south) -- (leaf02.north);
\node[branch,fill=white] (b1) at (0,1.6)
  {$\{1,4,7\}$\\[2pt]$y_1=x_1$};
\draw (root.south) -- node[right,pos=.56,fill=white,inner sep=1.5pt]
  {$x_0=1$} (b1.north);
\node[leaf] (leaf10) at (-1.35,0) {$1\mapsto 1$};
\draw (b1.south) -- (leaf10.north);
\node[leaf] (leaf11) at (0,0) {$4\mapsto 4$};
\draw (b1.south) -- (leaf11.north);
\node[leaf] (leaf12) at (1.35,0) {$7\mapsto 7$};
\draw (b1.south) -- (leaf12.north);
\node[branch,fill=black!10] (b2) at (5,1.6)
  {$\{2,5,8\}$\\[2pt]$y_1=x_1+1$};
\draw (root.south) -- node[above,sloped,pos=.56,fill=white,inner sep=1.5pt]
  {$x_0=2$} (b2.north);
\node[leaf] (leaf20) at (3.65,0) {$2\mapsto 5$};
\draw (b2.south) -- (leaf20.north);
\node[leaf] (leaf21) at (5,0) {$5\mapsto 8$};
\draw (b2.south) -- (leaf21.north);
\node[leaf] (leaf22) at (6.35,0) {$8\mapsto 2$};
\draw (b2.south) -- (leaf22.north);

\end{tikzpicture}
\caption{The permutation $f(x)=x+3[x\equiv2\pmod3]$ of
Example~\ref{ex:Z9}. We write $x=x_0+3x_1$ and $f(x)=y_0+3y_1$.
The branches group inputs by $x_0$. The equations for the digits are
in $\F_3$, and each leaf records a pair $x\mapsto f(x)$ in $\Z_9$.}
\label{fig:digit-tree-nine}
\end{figure}

Figure~\ref{fig:digit-tree-nine} shows the digit structure of the
non-affine permutation of $\Z_9$ in Example~\ref{ex:Z9}.
The first digit is unchanged. Within each residue class modulo
$3$, the second digit is translated, with the same multiplier $1$ in
every branch. The leaves give all nine values of $f$.

\section{Extremal permutations with several prime divisors}
\label{sec:mixed}

In this section, we study the extremal permutations of $\Z_n$ when $n$
has at least two distinct prime divisors. By
Theorem~\ref{thm:classification}, an extremal permutation of prime-power
order is of triangular type. We prove the same conclusion for a product
of two distinct odd primes in Theorem~\ref{thm:pq}, and we count these
permutations. For every other modulus with at least two distinct prime
divisors, we construct extremal permutations that are not of triangular
type. We begin with an exact criterion for $\varphi(n)$ good slopes. For
$z\in\Z_n$, we write $\gcd(z,n)$ for the greatest common divisor of $n$
and any integer representative of $z$. In particular, $\gcd(0,n)=n$.

\begin{lemma}\label{lem:gcd}
Let $g:\Z_n\to\Z_n$ be a map. Then, $x\mapsto g(x)+ux$ is a permutation of
$\Z_n$ for every $u\in\Z_n^\times$ if and only if
$\gcd(x-y,n)\neq\gcd(g(x)-g(y),n)$ for all distinct $x,y\in\Z_n$.
Consequently, a permutation $f$ of $\Z_n$ has $\varphi(n)$ good slopes if
and only if there is $a\in\Z_n^\times$ such that
$\gcd(x-y,n)\neq\gcd\bigl(f(x)-f(y)-a(x-y),n\bigr)$ for all distinct
$x,y\in\Z_n$.
\end{lemma}

\begin{proof}
We first show that, for $d,h\in\Z_n$, the equation $ud=h$ has a solution
$u\in\Z_n^\times$ if and only if $\gcd(d,n)=\gcd(h,n)$. Multiplication by a
unit does not change the greatest common divisor with $n$, and this proves
necessity. For sufficiency, let $p^e\,\|\,n$. If $p^e$ divides
$\gcd(d,n)=\gcd(h,n)$, then every unit $u_p$ of $\Z_{p^e}$ satisfies
$u_pd\equiv h\pmod{p^e}$. Otherwise, the largest power of $p$ that
divides $\gcd(d,n)$ is $p^v$ for some $v$ with $0\le v<e$. Then,
$d\equiv p^vd_0$ and $h\equiv p^vh_0\pmod{p^e}$ with $d_0$ and $h_0$
prime to $p$. Every integer $u_p$ with $u_p\equiv h_0d_0^{-1}\pmod{p^{e-v}}$
is prime to $p$ and satisfies $u_pd\equiv h\pmod{p^e}$. By the Chinese
remainder theorem, these elements $u_p$ combine into a unit $u$ of $\Z_n$
with $ud=h$.

The map $x\mapsto g(x)+ux$ fails to be injective exactly when
$u(x-y)=-(g(x)-g(y))$ for some distinct $x$ and $y$. Since
$\gcd(-z,n)=\gcd(z,n)$, the first paragraph shows that some unit $u$
satisfies this equation for a given pair if and only if
$\gcd(x-y,n)=\gcd(g(x)-g(y),n)$. This proves the first assertion. For the
second, let $a\in\Z_n^\times$ and put $g(x)=f(x)-ax$. Then,
$f(x)+mx=g(x)+(m+a)x$. Hence, $\Z_n^\times-a\subseteq D_n(f)$ holds if and
only if $g(x)+ux$ is a permutation for every unit $u$. By
Theorem~\ref{thm:universal}, $|D_n(f)|=\varphi(n)$ holds if and only if
$D_n(f)=\Z_n^\times-a$ for some unit $a$, and this is equivalent to
$\Z_n^\times-a\subseteq D_n(f)$.
\end{proof}

The next lemma is a counting principle for weights on the lines of
$\F_r^2$, where $r$ is a prime. For a map $w:\F_r^2\to\Z$, the
\emph{weight} of a line of $\F_r^2$ is the sum of the values of $w$ on
its points. We use the lemma in the proof of Theorem~\ref{thm:pq}, and
again in Section~\ref{sec:gap}.

\begin{lemma}\label{lem:uniform}
Let $r$ be a prime.
\begin{enumerate}
\item Let $w:\F_r^2\to\Z$ be a map and $L\in\Z$. Suppose that, for $r$
of the $r+1$ directions of $\F_r^2$, every line with that direction has
weight $L$. Then, $w$ is constant on every line with the
remaining direction.
\item Let $T$ be a finite set of size $L$, and let $A:\F_r\times T\to\F_r$
be a map. Suppose that, for every $a\in\F_r^\times$ and $c\in\F_r$,
exactly $L$ pairs $(s,t)$ satisfy $A(s,t)+as=c$. Then, for each
$z\in\F_r$, the number $|\{t\in T:A(s,t)=z\}|$ is independent of $s$.
\end{enumerate}
\end{lemma}

\begin{proof}
For (i), the plane $\F_r^2$ is the disjoint union of the $r$ lines with
any of the specified directions. Hence, its total weight is $rL$. Fix a
point $P$, and let $\ell$ be the line through $P$ with the remaining
direction. The other $r$ lines through $P$ have weight $L$ each. The
point $P$ lies on all of them, every point outside $\ell$ lies on
exactly one of them, and the other points of $\ell$ lie on none of
them. Summing their weights, we obtain
$$
rL=r\,w(P)+rL-\sum_{Q\in\ell}w(Q).
$$
Hence, $r\,w(P)=\sum_{Q\in\ell}w(Q)$. The right-hand side does not
depend on the choice of $P$ in $\ell$, and therefore $w$ is constant on
$\ell$.

For (ii), put $w(s,z)=|\{t\in T:A(s,t)=z\}|$. Every vertical line
$s=s_0$ has weight $|T|=L$. By hypothesis, every line $z+as=c$ with
$a\neq0$ also has weight $L$. These are all directions except the
horizontal one. By (i), $w$ is constant on every horizontal line, which
is the assertion.
\end{proof}

We can now classify the extremal permutations of $\Z_{pq}$.

\begin{theorem}\label{thm:pq}
Let $p,q$ be distinct odd primes. A permutation of
$\Z_{pq}\cong\F_p\times\F_q$ is extremal if and only if it has one of
the forms
$$
f(s,t)=(as+c,\,bt+\theta(s))
\quad\text{or}\quad
f(s,t)=(as+\eta(t),\,bt+d),
$$
where $a\in\F_p^\times$, $b\in\F_q^\times$, $c\in\F_p$, $d\in\F_q$,
and $\theta:\F_p\to\F_q$, $\eta:\F_q\to\F_p$ are arbitrary maps.
In particular, a permutation of $\Z_{pq}$ is extremal if and only if it
is of triangular type. In either form,
$$
D_{pq}(f)=\{(\alpha,\beta):\alpha\ne-a,\ \beta\ne-b\},
\quad \mu_{pq}(f)=\min(p,q)-1.
$$
\end{theorem}

\begin{proof}
We may assume $p<q$. Let $f$ be extremal. By
Theorem~\ref{thm:universal}, there is a unit $c_0$ with
$D_{pq}(f)=\Z_{pq}^\times-c_0$. Put $g(x)=f(x)-c_0x$, and write
$g(s,t)=(g_1(s,t),g_2(s,t))$. Then,
$g(s,t)+(\alpha s,\beta t)$ is a permutation for every
$\alpha\in\F_p^\times$ and $\beta\in\F_q^\times$.

Lemma~\ref{lem:gcd} gives three consequences. First, if both input
coordinates differ, then at least one output coordinate of $g$ agrees.
Second, for two inputs with the same first coordinate, equality of the
first output coordinates forces equality of the second output
coordinates. Third, for two inputs with the same second coordinate,
equality of the second output coordinates forces equality of the first
output coordinates. In each case, the excluded alternative would give
the same pattern of zero and nonzero coordinates for the input and
output differences. Their greatest common divisors with $pq$ would then
be equal.

For each $\alpha\ne0$, the first coordinate $g_1(s,t)+\alpha s$ of a
permutation takes each value exactly $q$ times. By
Lemma~\ref{lem:uniform}(ii), there are integers $c_z\ge0$, where
$z\in\F_p$, such that $|\{t\in\F_q:g_1(s,t)=z\}|=c_z$ for every
$s\in\F_p$. Moreover, $\sum_{z\in\F_p}c_z=q$.
If $g_1$ is constant, the second consequence shows that
$g_2(s,t)$ depends only on $s$. This is the first required form after
adding $c_0x$.

Suppose that $g_1$ is not constant, and put
$C=\{z\in\F_p:c_z>0\}$. Then, $|C|\ge2$. Since $q>p$, some
$z_0\in C$ satisfies $c_{z_0}\ge2$. For every $s$, the second
consequence defines a function $h_s:C\to\F_q$ by
$g_2(s,t)=h_s(g_1(s,t))$. Let $s\ne s'$ and $z\in C\setminus\{z_0\}$.
Choose $t'$ with $g_1(s',t')=z$. There are at least two choices of $t$
with $g_1(s,t)=z_0$. Hence, we can choose one different from $t'$. Both
input coordinates now differ, and the first output coordinates differ.
The first consequence gives $h_s(z_0)=h_{s'}(z)$.

Fix $z\in C\setminus\{z_0\}$. Given two values $s_1,s_2$, choose
$s'$ distinct from both. This is possible because $p\ge3$.
Then, $h_{s_1}(z_0)=h_{s'}(z)=h_{s_2}(z_0)$. Thus, all $h_s(z_0)$
have a common value. The relation in the preceding paragraph implies
that every $h_s(z)$ has this value for every $z\in C$. Hence, $g_2$ is
constant. The third consequence now shows that $g_1(s,t)$
depends only on $t$. This gives the second required form.

Conversely, each of the two forms is the sum of the affine map
$(s,t)\mapsto(as,bt)$ and a map that is strictly triangular with
respect to a chain $\Z_{pq}\supsetneq I\supsetneq\{0\}$. For the first
form, this map is $(s,t)\mapsto(c,\theta(s))$ and $I=\{0\}\times\F_q$.
For the second form, it is $(s,t)\mapsto(\eta(t),d)$ and
$I=\F_p\times\{0\}$. Proposition~\ref{prop:triangular} gives the
asserted sets of good slopes and the value of $\mu_{pq}$.
\end{proof}

\begin{corollary}\label{cor:pqcount}
Let $p,q$ be distinct odd primes. The number of extremal permutations
of $\Z_{pq}$ is
$$
pq(p-1)(q-1)\bigl(q^{p-1}+p^{q-1}-1\bigr).
$$
Exactly $pq(p-1)(q-1)$ of them are affine. The number of extremal
permutations satisfying $f(0)=0$ and $f(1)=1$ is
$$
(q-1)q^{p-2}+(p-1)p^{q-2}-1,
$$
and the identity is the only affine permutation among them.
\end{corollary}

\begin{proof}
In the first form of Theorem~\ref{thm:pq}, there are $p-1$ choices for
$a$, $q-1$ choices for $b$, $p$ choices for $c$, and $q^p$ choices for
$\theta$. These parameters are uniquely determined by the map.
The second form gives $q(p-1)(q-1)p^q$ maps. A map belongs to both
families exactly when both coordinates depend only on their respective
input coordinates. This means that $\theta$ and $\eta$ are constant,
and the intersection consists of the $pq(p-1)(q-1)$ affine permutations.
Adding the two counts and subtracting the intersection proves the first
formula.

We call a permutation \emph{normalised} if $f(0)=0$ and $f(1)=1$. In
the first form, $f(0)=0$ forces $c=0$ and $\theta(0)=0$.
The condition $f(1)=1$ then forces $a=1$ and $b=1-\theta(1)$.
There are $q-1$ choices for $\theta(1)\ne1$, and $q^{p-2}$ choices for
the remaining values of $\theta$. The second form contributes
$(p-1)p^{q-2}$ maps. The identity is the only normalised map in both
families and the only normalised affine permutation. This proves the
last formula.
\end{proof}

\begin{example}\label{ex:pqcount15}
By Corollary~\ref{cor:pqcount}, there are
$15\cdot8(5^2+3^4-1)=12600$ extremal permutations of $\Z_{15}$.
Of these, $120$ are affine and $12480$ are non-affine. Under the
normalisation $f(0)=0$, $f(1)=1$, the two triangular families contain
$4\cdot5=20$ and $2\cdot3^3=54$ maps. They meet only in the identity.
Thus, there are $73$ normalised extremal permutations, of which $72$
are non-affine.
\end{example}

Strictly triangular maps satisfy the criterion of Lemma~\ref{lem:gcd} for
a structural reason. The ideal generated by $z\in\Z_n$ is $\gcd(z,n)\Z_n$.
Let $g$ be strictly triangular with respect to a chain of ideals $I_k$.
If $x-y\in I_k\setminus I_{k+1}$, then $g(x)-g(y)\in I_{k+1}$. Hence, the two
differences generate different ideals, and their greatest common divisors
with $n$ are different.

The next lemma gives a simple way to show that a map is not strictly
triangular with respect to any chain.

\begin{lemma}\label{lem:obstruction}
Let $g:\Z_n\to\Z_n$ be a map, and let $P$ be a set of pairs of elements of
$\Z_n$. Let $I$ be the ideal generated by the differences $x-y$, and let
$J$ be the ideal generated by the differences $g(x)-g(y)$, where $(x,y)$
runs over $P$. If $I\neq\{0\}$ and $I\subseteq J$, then $g$ is not
strictly triangular with respect to any chain of ideals.
\end{lemma}

\begin{proof}
Suppose that $g$ is strictly triangular with respect to
$\Z_n=I_0\supsetneq\cdots\supsetneq I_r=\{0\}$. Since $I\subseteq I_0$ and
$I\not\subseteq I_r$, there is a largest index $k<r$ with $I\subseteq I_k$.
For every $(x,y)\in P$, we have $x-y\in I_k$, and hence
$g(x)-g(y)\in I_{k+1}$. Thus, $J\subseteq I_{k+1}$, and
$I\subseteq J\subseteq I_{k+1}$. This contradicts the choice of $k$.
\end{proof}

Adding a constant to $g$ does not change the differences $g(x)-g(y)$. Hence,
a permutation $f$ is not of triangular type if, for every $a\in\Z_n$, the
map $x\mapsto f(x)-ax$ satisfies the hypothesis of
Lemma~\ref{lem:obstruction} for some set of pairs.

We now construct extremal permutations that are not of triangular type.
Let $n=n_1n_2$ with $\gcd(n_1,n_2)=1$. We identify $\Z_n$ with
$\Z_{n_1}\times\Z_{n_2}$. The elements $(1,0)$ and $(0,1)$ are idempotents
whose sum is $1$. Hence, the ideal generated by an element $(z_1,z_2)$
contains $(z_1,0)$ and $(0,z_2)$. Moreover, for a divisor $d_1$ of
$n_1$, the ideal generated by $(d_1,0)$ and $(0,1)$ is
$d_1\Z_{n_1}\times\Z_{n_2}$.

\begin{proposition}\label{prop:switch}
Let $n=n_1n_2$ with $\gcd(n_1,n_2)=1$ and $n_2\ge2$, and let $d$ be a
divisor of $n_1$ with $1<d<n_1$. On $\Z_n=\Z_{n_1}\times\Z_{n_2}$, define
$$
f(x,y)=\bigl(x+d\,[x\equiv0\ (\mathrm{mod}\ d)]\,[y=1],\ y+[x=1+d]\bigr).
$$
Then, $f$ is a permutation of $\Z_n$ with $D_n(f)=\Z_n^\times-1$ and
$\mu_n(f)=p_0(n)-1$, and $f$ is not of triangular type.
\end{proposition}

\begin{proof}
Put $g(x,y)=f(x,y)-(x,y)=\bigl(d\,[x\equiv0\ (\mathrm{mod}\ d)]\,[y=1],\
[x=1+d]\bigr)$, and let $(v_1,v_2)$ be a unit of $\Z_n$. Then, $v_1$ and
$v_2$ are units of $\Z_{n_1}$ and $\Z_{n_2}$. We show that the map
$$
(x,y)\longmapsto(X,Y)=g(x,y)+(v_1x,v_2y)
=\bigl(v_1x+d\,[x\equiv0\ (\mathrm{mod}\ d)]\,[y=1],\ v_2y+[x=1+d]\bigr)
$$
is injective. Since $X\equiv v_1x\pmod d$ and $v_1$ is a unit modulo $d$,
the value $X$ determines $x$ modulo $d$. If $x\equiv0\pmod d$, then
$x\neq1+d$, because $d>1$. Hence, $Y=v_2y$ determines $y$, and then $X$
determines $x$. If $x\equiv1\pmod d$, then $X=v_1x$ determines $x$, and
then $Y$ determines $y$. For all other $x$, we have $(X,Y)=(v_1x,v_2y)$.
Thus, the map is injective, and therefore a permutation. Consequently,
$f(z)+mz=g(z)+(1+m)z$ is a permutation whenever $1+m$ is a unit. In
particular, $f$ is a permutation and $\Z_n^\times-1\subseteq D_n(f)$.
Theorem~\ref{thm:universal} gives $D_n(f)=\Z_n^\times-1$ and
$\mu_n(f)=p_0(n)-1$.

Let $a=(a_1,a_2)\in\Z_n$, and put $h(z)=f(z)-az$. Note that
$1+d\not\equiv0,1\pmod{n_1}$, since $d\le n_1/2$ and $n_1\ge4$. We use the
pairs $\bigl((0,1),(0,0)\bigr)$ and $\bigl((1+d,0),(1,0)\bigr)$. Their
differences are $(0,1)$ and $(d,0)$, which generate the ideal
$I=d\Z_{n_1}\times\Z_{n_2}\neq\{0\}$. We have $f(0,1)=(d,1)$, $f(0,0)=(0,0)$,
$f(1+d,0)=(1+d,1)$ and $f(1,0)=(1,0)$. Hence,
$$
h(0,1)-h(0,0)=(d,\,1-a_2),\quad h(1+d,0)-h(1,0)=\bigl((1-a_1)d,\,1\bigr).
$$
The ideal $J$ generated by these two elements contains $(d,0)$ and
$(0,1)$. Hence, $J\supseteq I$. By Lemma~\ref{lem:obstruction}, $h$ is not
strictly triangular with respect to any chain. As $a$ was arbitrary, $f$
is not of triangular type.
\end{proof}

For $n=2q$ with $q$ an odd prime, no factor $n_1$ as in
Proposition~\ref{prop:switch} exists. In this case, we use the fact that
every unit of $\Z_{2q}$ acts trivially on the factor $\Z_2$.

\begin{proposition}\label{prop:2q}
Let $q$ be an odd prime and $n=2q$. On $\Z_n=\Z_2\times\Z_q$, define
$f(s,t)=\bigl(s',\,t+[s'=1]\bigr)$, where $s'=s+[t=1]$. Then, $f$ is a
permutation of $\Z_n$ with $D_n(f)=\Z_n^\times-1$, and $f$ is not of
triangular type.
\end{proposition}

\begin{proof}
Put $g(s,t)=f(s,t)-(s,t)=\bigl([t=1],\,[s'=1]\bigr)$. The units of $\Z_n$
are the elements $(1,w)$ with $w\in\Z_q^\times$, and
$g(s,t)+(1,w)(s,t)=\bigl(s',\,wt+[s'=1]\bigr)$. The output determines $s'$,
then it determines $t$ through its second coordinate, and finally
$s=s'-[t=1]$. Thus, $g(z)+uz$ is a permutation for every unit $u$. As in
the proof of Proposition~\ref{prop:switch}, it follows that $f$ is a
permutation with $D_n(f)=\Z_n^\times-1$.

Let $a=(a_1,a_2)\in\Z_n$, and put $h(z)=f(z)-az$. We have $f(0,0)=(0,0)$,
$f(0,1)=(1,2)$ and $f(1,0)=(1,1)$. For the pairs
$\bigl((0,1),(0,0)\bigr)$ and $\bigl((1,0),(0,0)\bigr)$, we obtain
$$
h(0,1)-h(0,0)=(1,\,2-a_2),\quad h(1,0)-h(0,0)=(1-a_1,\,1).
$$
The ideal generated by these two elements contains $(1,0)$ and $(0,1)$.
Hence, it is all of $\Z_n$. Lemma~\ref{lem:obstruction} applies with $I=\Z_n$,
and $f$ is not of triangular type.
\end{proof}

For $q=3$, the permutation of Proposition~\ref{prop:2q} is the $4$-cycle
$(1\ 4\ 5\ 3)$ of $\Z_6$, and its good slopes are $0$ and $4$. The
transposition $(1\ 4)$ of Example~\ref{ex:coprime} has the same good
slopes, but it is of triangular type.

\begin{theorem}\label{thm:nontriangular}
Every extremal permutation of $\Z_n$ is of triangular type if and only
if $n$ is a prime power or a product of two distinct odd primes.
For every other $n$, there is an extremal permutation that is not of
triangular type.
\end{theorem}

\begin{proof}
The affirmative cases follow from Theorems~\ref{thm:classification}
and~\ref{thm:pq}. Suppose that $n$ is neither a prime power nor a product
of two distinct odd primes. Then, $n$ has at least two distinct prime
divisors. If $p^2\mid n$ and $p^e\,\|\,n$, apply
Proposition~\ref{prop:switch} with $n_1=p^e$, $n_2=n/p^e$ and $d=p$.
Here, $n_2\ge2$ because $n$ is not a prime power.

If $n$ is squarefree with at least three prime divisors, choose two of
them, say $p$ and $\ell$. Proposition~\ref{prop:switch} applies with
$n_1=p\ell$, $n_2=n/(p\ell)$ and $d=p$.
The remaining possibility is $n=2q$ with $q$ an odd prime.
Proposition~\ref{prop:2q} applies in this situation.
\end{proof}

\begin{example}\label{ex:45}
Let $n=45$, $n_1=9$, $n_2=5$ and $d=3$ in Proposition~\ref{prop:switch}.
The idempotents $(1,0)$ and $(0,1)$ of $\Z_9\times\Z_5$ correspond to $10$
and $36$ in $\Z_{45}$. For $z\in\Z_{45}$, we have $x=z\bmod 9$ and
$y=z\bmod 5$. Hence, the condition $x\equiv0\pmod 3$ together with $y=1$
means that $z\equiv6\pmod{15}$, and $x=4$ means that $z\equiv4\pmod 9$.
Thus, the permutation of Proposition~\ref{prop:switch} is
$$
f(z)=z+30\,[z\equiv6\ (\mathrm{mod}\ 15)]+36\,[z\equiv4\ (\mathrm{mod}\ 9)].
$$
In cycle notation, $f=(6\ 36\ 21)(4\ 40\ 31\ 22\ 13)$, and all other
elements are fixed. Its $24$ good slopes are the elements of $\Z_{45}$
that are congruent to $0$ or $1$ modulo $3$ and not congruent to $4$
modulo $5$.
The set $\{0,1\}$ is admissible, and $\mu_{45}(f)=2$. The obstruction of
the proof can be checked directly in $\Z_{45}$. For every $a\in\Z_{45}$,
the map $h(z)=f(z)-az$ satisfies $h(36)-h(0)=21-36a$ and
$h(40)-h(10)=21-30a$. Multiplying by the idempotents $10$ and $36$, we
obtain $10(21-36a)=30$ and
$36(21-30a)=36$ in $\Z_{45}$. Hence, the ideal generated by the two
differences of $h$ contains $30$ and $36$, and thus contains $3\Z_{45}$.
This is also the ideal generated by the input differences $36$ and $30$.
Note that $|D_{45}(f)|=24>45/2$. Moreover, $f$ is not
congruence-preserving, since $0\equiv36\pmod 9$ while $f(0)=0$ and
$f(36)=21$.
\end{example}

\section{Congruence-preserving permutations}\label{sec:half}

In this section, we classify the congruence-preserving permutations with
more than $n/2$ good slopes. When $\varphi(n)>n/2$, this threshold lies
below the maximum $\varphi(n)$. For some moduli, the classification
therefore includes permutations that are not extremal. We first
describe congruence-preserving maps through their components modulo prime
powers and through their digit maps. We then prove the classification in
Theorem~\ref{thm:half}, count the permutations of prime-power order, and
treat polynomial maps.

\begin{lemma}\label{lem:CRT}
Let $f$ be a congruence-preserving map of $\Z_n$. For every $p^e\,\|\,n$,
there is a congruence-preserving map $f_p$ of $\Z_{p^e}$ with
$f(x)\equiv f_p(x\bmod p^e)\pmod{p^e}$ for all $x\in\Z_n$. The map $f$ is a
permutation if and only if every $f_p$ is a permutation. In this case,
$D_n(f)$ consists of the $m\in\Z_n$ with $m\bmod p^e\in D_{p^e}(f_p)$ for
every $p^e\,\|\,n$.
\end{lemma}

\begin{proof}
Since $f$ preserves congruences modulo $p^e$, the residue of $f(x)$ modulo
$p^e$ depends only on $x$ modulo $p^e$. This defines $f_p$, and $f_p$ is
congruence-preserving because every divisor of $p^e$ divides $n$. In the
coordinates of the Chinese remainder theorem, $f$ acts componentwise as
$(f_p)_{p\mid n}$. A componentwise map is a permutation if and only if each
component is a permutation. For $m\in\Z_n$, the map $x\mapsto f(x)+mx$ is
again congruence-preserving, and its components are the maps
$x\mapsto f_p(x)+(m\bmod p^e)x$. This gives the description of $D_n(f)$.
\end{proof}

We call the maps $f_p$ the \emph{components} of $f$. Now let $p$ be a
prime, let $e\ge1$, and let $g$ be a congruence-preserving map of
$\Z_{p^e}$. We write $x$ and $g(x)$ in base $p$ with digits $x_j$ and
$y_j$, where $0\le j<e$, as in Section~\ref{sec:primepowers}. Since
$g(x)$ modulo $p^{j+1}$ depends only on $x$ modulo $p^{j+1}$, the digit
$y_j$ depends only on $x_0,\dots,x_j$. Hence, for every $j$ and every
$u=(x_0,\dots,x_{j-1})\in\F_p^{\,j}$, there is a map
$\sigma_{j,u}:\F_p\to\F_p$ with $y_j=\sigma_{j,u}(x_j)$. For $j=0$, the
tuple $u$ is empty. We call the maps $\sigma_{j,u}$ the \emph{digit maps}
of $g$. There are $1+p+\dots+p^{e-1}$ of them.

\begin{lemma}\label{lem:digits}
Let $p$ be a prime, let $e\ge1$, and let $g$ be a congruence-preserving
map of $\Z_{p^e}$ with digit maps $\sigma_{j,u}$.
\begin{enumerate}
\item The map $g$ is a permutation if and only if every $\sigma_{j,u}$ is a
permutation of $\F_p$.
\item Let $g$ be a permutation, and let $T(g)$ be the intersection of the
sets $D_p(\sigma_{j,u})$ over all digit maps. Then,
$D_{p^e}(g)=\{m\in\Z_{p^e}:m\bmod p\in T(g)\}$. Moreover,
$|D_{p^e}(g)|=p^{e-1}|T(g)|$ and $\mu_{p^e}(g)=|T(g)|$.
\end{enumerate}
\end{lemma}

\begin{proof}
If every $\sigma_{j,u}$ is a permutation, the input digits can be
recovered one at a time from the output digits. Hence, $g$ is injective,
and therefore a permutation. Conversely, let $g$ be a permutation. The map
induced by $g$ on $\Z_{p^{j+1}}$ is surjective, and hence bijective.
Suppose that $\sigma_{j,u}(t)=\sigma_{j,u}(t')$ with $t\neq t'$. The two
residues modulo $p^{j+1}$ with lower digits $u$ and $j$th digits $t$ and
$t'$ then have the same image modulo $p^{j+1}$, which is impossible.

Let $m\in\Z_{p^e}$, and put $\bar m=m\bmod p$. The map
$x\mapsto g(x)+mx$ is again congruence-preserving, and we compute its
digit maps. Fix $j$ and $u$, and put $U=\sum_{i<j}x_ip^i$. Then,
$x\equiv U+p^jx_j\pmod{p^{j+1}}$. The lower digits $y_0,\dots,y_{j-1}$ of
$g(x)$ depend only on $u$. Hence,
$g(x)\equiv G_u+p^j\sigma_{j,u}(x_j)\pmod{p^{j+1}}$ with
$G_u=\sum_{i<j}y_ip^i$. Moreover, $mx\equiv mU+p^j\bar mx_j\pmod{p^{j+1}}$.
Write $G_u+mU=A+p^jB$ with integers $0\le A<p^j$ and $B$. Then,
$$
g(x)+mx\equiv A+p^j\bigl(\sigma_{j,u}(x_j)+\bar mx_j+B\bigr)
\pmod{p^{j+1}}.
$$
Thus, the digit map of $g(x)+mx$ for $j$ and $u$ is
$t\mapsto\sigma_{j,u}(t)+\bar mt+B$, where $B$ depends only on $u$ and $m$.
By (i), $g(x)+mx$ is a permutation if and only if
$t\mapsto\sigma_{j,u}(t)+\bar mt$ is a permutation of $\F_p$ for all $j$
and $u$, that is, if and only if $\bar m\in T(g)$. This proves the
description of $D_{p^e}(g)$. The set $T(g)$ contains $0$ by (i), and the
remaining assertions follow from Lemma~\ref{lem:admissible}(ii).
\end{proof}

The following proposition explains why congruence preservation is natural
for permutations of prime-power order with many good slopes.

\begin{proposition}\label{prop:fullclass}
Let $p$ be a prime, let $e\ge1$, and let $f$ be a permutation of
$\Z_{p^e}$. The following are equivalent.
\begin{enumerate}
\item $f$ is congruence-preserving.
\item $p\Z_{p^e}\subseteq D_{p^e}(f)$.
\item $D_{p^e}(f)$ contains a full residue class modulo $p$.
\end{enumerate}
Consequently, if $f$ is not congruence-preserving, then
$|D_{p^e}(f)|\le(p-1)(p^{e-1}-1)$.
\end{proposition}

\begin{proof}
If $f$ is congruence-preserving, its digit maps are permutations by
Lemma~\ref{lem:digits}(i). Hence, $0\in T(f)$, and
Lemma~\ref{lem:digits}(ii) gives (ii). Clearly, (ii) implies (iii).

Now assume (iii), say $c+p\Z_{p^e}\subseteq D_{p^e}(f)$. By
Lemma~\ref{lem:normal}, the map $g(x)=f(x)+cx$ is a permutation with
$p\Z_{p^e}\subseteq D_{p^e}(g)$. Let $h$ be the inverse of $g$. We claim
that $h$ is congruence-preserving. Otherwise, there are $y,z\in\Z_{p^e}$
and $k\ge1$ with $y\equiv z\pmod{p^k}$ and $h(y)\not\equiv h(z)\pmod{p^k}$.
Then, $y\neq z$, and $v=\nu(h(y)-h(z))$ and $w=\nu(y-z)$ satisfy
$v<k\le w$. Write $h(y)-h(z)=p^vc_0$ and $y-z=p^wd_0$, where $c_0$ is prime
to $p$. Then, $m=-p^{w-v}d_0c_0^{-1}$ lies in $p\Z_{p^e}$, and it satisfies
$m\bigl(h(y)-h(z)\bigr)=-(y-z)$, that is, $y+mh(y)=z+mh(z)$. Since
$g(h(y))=y$ and $g(h(z))=z$, the map $x\mapsto g(x)+mx$ takes the same
value at the distinct elements $h(y)$ and $h(z)$. This contradicts
$m\in D_{p^e}(g)$, and the claim follows.

For every $k$, let $h_k$ be the map induced by $h$ on $\Z_{p^k}$. It is
surjective, and hence bijective. If $x\equiv x'\pmod{p^k}$, then
$h_k\bigl(g(x)\bmod p^k\bigr)=x\bmod p^k=h_k\bigl(g(x')\bmod p^k\bigr)$, and
therefore $g(x)\equiv g(x')\pmod{p^k}$. Thus, $g$ is
congruence-preserving, and so is $f(x)=g(x)-cx$. This proves (i).

Finally, suppose that $f$ is not congruence-preserving. Since (iii)
implies (i), no residue class modulo $p$ is contained in $D_{p^e}(f)$. By
Theorem~\ref{thm:universal}, $D_{p^e}(f)$ misses one residue class
completely, and it misses at least one element of every other class.
Each class has $p^{e-1}$ elements, and the bound follows.
\end{proof}

\begin{corollary}\label{cor:Z9half}
A permutation of $\Z_9$ has more than $9/2$ good slopes if and only if it
is one of the $162$ permutations of Example~\ref{ex:Z9}.
\end{corollary}

\begin{proof}
By Proposition~\ref{prop:fullclass}, a permutation of $\Z_9$ that is not
congruence-preserving has at most $2\cdot2=4$ good slopes. For a
congruence-preserving permutation, Lemma~\ref{lem:digits} shows that the
number of good slopes is a multiple of $3$, and Theorem~\ref{thm:universal}
bounds it by $6$. Hence, a number of good slopes larger than $9/2$ must be
$6=\varphi(9)$. Conversely, the permutations of Example~\ref{ex:Z9} have
six good slopes.
\end{proof}

We now treat congruence-preserving permutations of arbitrary order. Let
$f$ be a congruence-preserving permutation of $\Z_n$ with components $f_p$.
The digit maps of $f_p$ are permutations of $\F_p$ by
Lemma~\ref{lem:digits}(i). If all of them are affine, we write $A_p$ for
the set of their multipliers. Thus, $A_p$ is a nonempty subset of
$\F_p^\times$.

\begin{theorem}\label{thm:half}
Let $f$ be a congruence-preserving permutation of $\Z_n$ with components
$f_p$. Then, $|D_n(f)|>n/2$ if and only if the following hold.
\begin{enumerate}
\item $n$ is odd.
\item Every digit map of every component $f_p$ is affine.
\item The sets $A_p$ of multipliers satisfy
$\prod_{p\mid n}\bigl(1-|A_p|/p\bigr)>1/2$.
\end{enumerate}
In this case, $D_n(f)$ consists of the $m\in\Z_n$ with
$m\bmod p\notin-A_p$ for every prime $p\mid n$, and
$$
|D_n(f)|=n\prod_{p\mid n}\Bigl(1-\frac{|A_p|}{p}\Bigr),\quad
\mu_n(f)=\min_{p\mid n}\bigl(p-|A_p|\bigr).
$$
\end{theorem}

\begin{proof}
Put $T_p=T(f_p)$, with the notation of Lemma~\ref{lem:digits}. By
Lemmas~\ref{lem:CRT} and~\ref{lem:digits}, the set $D_n(f)$ consists of
the $m\in\Z_n$ with $m\bmod p\in T_p$ for every $p\mid n$. Hence,
Lemma~\ref{lem:admissible}(ii) gives
$$
|D_n(f)|=n\prod_{p\mid n}\frac{|T_p|}{p}\quad\text{and}\quad
\mu_n(f)=\min_{p\mid n}|T_p|.
$$
Suppose that $|D_n(f)|>n/2$. Each factor $|T_p|/p$ is at most $1$.
Hence, every factor exceeds $1/2$. For $p=2$, this is impossible. Indeed,
Theorem~\ref{thm:universal} gives $2^{e-1}|T_2|=|D_{2^e}(f_2)|\le2^{e-1}$,
where $2^e\,\|\,n$. Hence, $n$ is odd. Every digit map $\sigma$ of $f_p$
satisfies $D_p(\sigma)\supseteq T_p$. Hence, it has more than $p/2$ good
slopes. By Theorem~\ref{thm:fields}(ii), $\sigma$ is affine. An affine
permutation $t\mapsto\alpha t+\beta$ of $\F_p$ has the good slopes
$\F_p\setminus\{-\alpha\}$. Hence, $T_p=\F_p\setminus(-A_p)$ and
$|T_p|=p-|A_p|$, which gives (iii) and the stated formulas. Conversely, if
all digit maps are affine, the same computation gives
$T_p=\F_p\setminus(-A_p)$, and (iii) gives $|D_n(f)|>n/2$.
\end{proof}

For an odd prime power $n=p^e$, condition (iii) reads $|A_p|<p/2$, and
Theorem~\ref{thm:classification} is the case $|A_p|=1$. Moreover, for
odd $n$, every choice of affine digit maps with nonzero multipliers
defines a congruence-preserving permutation of each $\Z_{p^e}$ by
Lemma~\ref{lem:digits}(i), and hence of $\Z_n$ by the Chinese remainder
theorem. Thus, condition (iii) is the only restriction on these
choices.
An affine permutation $x\mapsto ax+b$ has all multipliers equal to
$a\bmod p$. Hence, if $|A_p|\ge2$ for some $p$, then $f$ is not affine.

\begin{corollary}\label{cor:counthalf}
Let $p$ be an odd prime, let $e\ge1$, put $J=1+p+\dots+p^{e-1}$, and let
$s$ be an integer with $1\le s\le(p+1)/2$. The number of
congruence-preserving permutations of $\Z_{p^e}$ with $p^{e-1}(p-s)$ good
slopes is
$$
\binom{p-1}{s}\,p^J\sum_{k=0}^{s}(-1)^k\binom sk(s-k)^J.
$$
Each such permutation $f$ satisfies $\mu_{p^e}(f)=p-s$.
\end{corollary}

\begin{proof}
By Lemma~\ref{lem:digits}, such a permutation $f$ satisfies
$|T(f)|=p-s\ge(p-1)/2$. Hence, every digit map has more than $(p-3)/2$
good slopes, and Theorem~\ref{thm:fields}(i) shows that it is affine.
The set of good slopes of an affine digit map omits only the negative of
its multiplier. Hence, $|T(f)|=p-s$ exactly when there are $s$ distinct
multipliers.
Conversely, any such collection of affine digit maps defines a
permutation with $p^{e-1}(p-s)$ good slopes. There are $J$ digit maps.
We choose the set of multipliers in $\binom{p-1}{s}$ ways, and then a
surjective assignment of these multipliers to the $J$ digit maps. By inclusion
and exclusion, the
number of such assignments is the alternating sum in the statement. The
constants give $p^J$ further independent choices. The digit maps are
determined by $f$. Hence, no permutation is counted twice. The value of
$\mu_{p^e}(f)$ comes from Lemma~\ref{lem:digits}.
\end{proof}

\begin{example}\label{ex:Z25}
Let $n=25$. Then, $J=6$, and for $s=1$, Corollary~\ref{cor:counthalf} gives
$4\cdot5^6=62500$ permutations with $20$ good slopes, in accordance with
Corollary~\ref{cor:count}. For $s=2$, it gives
$6\cdot5^6\cdot(2^6-2)=5812500$ congruence-preserving permutations $f$
with $15$ good slopes, and each of them satisfies $\mu_{25}(f)=3$.

Now let $p=3$ and $e\ge1$. Every permutation of $\F_3$ is affine. Hence,
by Lemma~\ref{lem:digits}, a congruence-preserving permutation of
$\Z_{3^e}$ has $3^{e-1}(3-s)$ good slopes, where $s\in\{1,2\}$ is the
number of distinct multipliers of its digit maps. Thus, if it has more
than $3^e/2$ good slopes, then $s=1$, and it has $\varphi(3^e)$ good
slopes.
\end{example}

For polynomial maps, the digit maps can be computed explicitly. For
$F\in\Z[X]$ and a prime $p$, let $\overline F_p:\F_p\to\F_p$ be the map
induced by $F$, and let $F'$ be the formal derivative of $F$.

\begin{proposition}\label{prop:poly}
Let $F\in\Z[X]$ induce the map $f$ of $\Z_n$, and let $p^{e_p}\,\|\,n$ for
every prime $p\mid n$.
\begin{enumerate}
\item The map $f$ is a permutation if and only if, for every prime
$p\mid n$, the map $\overline F_p$ is a permutation of $\F_p$ and, when
$e_p\ge2$, $F'(r)\not\equiv0\pmod p$ for all $r\in\Z$.
\item Let $f$ be a permutation. For every prime $p\mid n$, put
$$
T_p=\begin{cases}
D_p(\overline F_p)&\text{if }e_p=1,\\
D_p(\overline F_p)\setminus\{-F'(r)\bmod p:\ r\in\Z\}&\text{if }e_p\ge2.
\end{cases}
$$
Then, $D_n(f)$ consists of the $m\in\Z_n$ with $m\bmod p\in T_p$ for every
prime $p\mid n$.
\end{enumerate}
\end{proposition}

\begin{proof}
The map $f$ is congruence-preserving, and its component $f_p$ is the map
of $\Z_{p^e}$ induced by $F$, where $e=e_p$. Its digit map for $j=0$ is
$\overline F_p$. Let $e\ge2$, let $1\le j<e$, and let
$u=(x_0,\dots,x_{j-1})$ with $x_0=r$. As in the proof of
Lemma~\ref{lem:digits}, we have $x\equiv U+p^jx_j\pmod{p^{j+1}}$. Taylor
expansion with integral coefficients and $2j\ge j+1$ give
$F(x)\equiv F(U)+p^jx_jF'(U)\pmod{p^{j+1}}$. Moreover,
$F'(U)\equiv F'(r)\pmod p$. As in the proof of Lemma~\ref{lem:digits}, the
digit map $\sigma_{j,u}$ is therefore $t\mapsto F'(r)t+c_u$, where $c_u$
depends only on $u$. Every $r\in\F_p$ occurs as the first entry of such a
tuple $u$. By Lemmas~\ref{lem:CRT} and~\ref{lem:digits}(i), $f$ is a
permutation if and only if every $\overline F_p$ is a permutation and, when
$e_p\ge2$, every map $t\mapsto F'(r)t+c_u$ is a permutation. The latter
means $F'(r)\not\equiv0\pmod p$, and this proves (i).

For (ii), let $m\in\Z_n$. The map $x\mapsto f(x)+mx$ is induced by
$F(X)+mX$, whose reduction modulo $p$ induces $x\mapsto\overline F_p(x)+mx$
and whose derivative is $F'+m$. By (i), $f(x)+mx$ is a permutation if and
only if, for every $p\mid n$, $m\bmod p\in D_p(\overline F_p)$ and, when
$e_p\ge2$, $m\not\equiv-F'(r)\pmod p$ for all $r$. This is the stated
description.
\end{proof}

Part (i) recovers the well-known criterion for polynomial permutations of
$\Z_{p^e}$. For $p=2$ and $e\ge2$, Rivest expressed this criterion
through the parities of the coefficients of $F$ \cite{Rivest}.
The proof shows that the digit maps of $f_p$ are
$\overline F_p$ and the affine maps $t\mapsto F'(r)t+c_u$. Combined with
Theorem~\ref{thm:half}, this gives the following criterion. A polynomial
permutation has more than $n/2$ good slopes if and only if $n$ is odd,
every $\overline F_p$ is an affine map $x\mapsto a_px+b_p$, and
$\prod_{p\mid n}(1-|A_p|/p)>1/2$, where
$$
A_p=\begin{cases}
\{a_p\}&\text{if }e_p=1,\\
\{a_p\}\cup\{F'(r)\bmod p:\ r\in\Z\}&\text{if }e_p\ge2.
\end{cases}
$$

\begin{corollary}\label{cor:polyext}
There is a non-affine polynomial permutation $f$ of $\Z_n$ with
$|D_n(f)|=\varphi(n)$ if and only if $8\mid n$ or $p^2\mid n$ for some odd
prime $p$.
\end{corollary}

\begin{proof}
Let $f$ be a polynomial permutation with $|D_n(f)|=\varphi(n)$. By
Lemma~\ref{lem:CRT}, $|D_n(f)|$ is the product of the numbers
$|D_{p^e}(f_p)|$ over all $p^e\,\|\,n$. Each of them is at most
$\varphi(p^e)$ by Theorem~\ref{thm:universal}. Hence, every component $f_p$
has $\varphi(p^e)$ good slopes. Suppose that $8\nmid n$ and $p^2\nmid n$
for every odd prime $p$. Then, every $p^e\,\|\,n$ satisfies $e=1$ or
$(p,e)=(2,2)$. By Corollary~\ref{cor:count}, every component is then an
affine map $x\mapsto a_px+b_p$. Choosing $a,b\in\Z_n$ with $a\equiv a_p$
and $b\equiv b_p\pmod{p^e}$ for all $p$, we obtain $f(x)=ax+b$.

Conversely, let $p^e\,\|\,n$, where $p$ is odd and $e\ge2$, or $p=2$ and
$e\ge3$. Choose $c\in\Z$ with $c\equiv p\pmod{p^e}$ and $c\equiv0$ modulo
every other prime power exactly dividing $n$, and put $F(X)=X+cX^2$. On the
component $\Z_{p^e}$, $F$ induces $x+px^2$, which is a non-affine
permutation with $\varphi(p^e)$ good slopes by Example~\ref{ex:primepower}.
On every other component, it induces the identity. Hence, $F$ induces a
non-affine permutation with $\varphi(n)$ good slopes.
\end{proof}

In particular, if neither $8$ nor the square of an odd prime divides $n$,
then the non-affine permutations with $\varphi(n)$ good slopes in
Example~\ref{ex:coprime} are not induced by polynomials.

\begin{example}\label{ex:2xx5}
Let $n=25$, let $F(X)=2X-X^5$, and let $f$ be the map of $\Z_{25}$
induced by $F$. Since $x^5=x$ in $\F_5$, we have
$\overline F_5(x)=x$, and $F'(X)=2-5X^4\equiv2\pmod5$. By
Proposition~\ref{prop:poly}, $f$ is a permutation with
$T_5=\F_5\setminus\{-1,-2\}=\{0,1,2\}$. Thus,
$$
D_{25}(f)=\{m\in\Z_{25}:m\bmod5\in\{0,1,2\}\},\quad
|D_{25}(f)|=15>\frac{25}{2},\quad \mu_{25}(f)=3.
$$
It is not affine, since $f(2)-2f(1)+f(0)=-30\equiv20\pmod{25}$.

The same idea works for every odd $n$ whose largest prime divisor $q$
satisfies $q^2\mid n$. Let $q^{e}\,\|\,n$, choose $c\in\Z$ with
$c\equiv1\pmod{q^{e}}$ and $c\equiv0$ modulo every other prime power
exactly dividing $n$, and let $f$ be the map of $\Z_n$ induced by
$F(X)=X+c(X-X^q)$. Then, $\overline F_p(x)=x$ for every prime $p\mid n$.
Moreover, $F'\equiv2\pmod q$, and $F'\equiv1\pmod p$ for $p\neq q$. By
Proposition~\ref{prop:poly}, we obtain
$T_q=\F_q\setminus\{-1,-2\}$ and $T_p=\F_p\setminus\{-1\}$ for $p\neq q$.
Hence, $f$ is a permutation with $|D_n(f)|=\varphi(n)(q-2)/(q-1)$. It is
not affine, because affine permutations have $\varphi(n)$ good slopes. For
instance, $n=75$ and $c=51$ give a permutation of $\Z_{75}$ with $30$ good
slopes, whereas $\varphi(75)=40$.
\end{example}

\begin{example}\label{ex:notCP}
The permutation $f(x)=x+10\,[x\equiv1\pmod5]$ of $\Z_{15}$ from
Example~\ref{ex:coprime} has $8>15/2$ good slopes. However,
$0\equiv6\pmod 3$, while $f(0)=0$ and $f(6)=1$ are not congruent modulo
$3$. Thus, $f$ is not congruence-preserving,
and Theorem~\ref{thm:half} does not apply to it.
\end{example}

For squarefree $n$, Theorem~\ref{thm:half} shows that every
congruence-preserving permutation with more than $n/2$ good slopes is
affine. Indeed, each component then has a single digit map, which is
affine, and the components combine into one affine map. Hence, for
squarefree $n$, all non-affine permutations with more than $n/2$ good slopes,
such as the one in Example~\ref{ex:notCP}, are not congruence-preserving.

\section{The gap below the maximum}\label{sec:gap}

Let $q$ be the largest prime divisor of $n$, and put
$G(n)=\varphi(n)(q-2)/(q-1)$. We study whether a non-extremal
permutation can have more than $G(n)$ good slopes. We prove the bound
for congruence-preserving permutations and for several families of
moduli. We also give a reduction through preserved congruences and
construct permutations attaining $G(n)$ for every odd composite $n$.

For congruence-preserving permutations, the digit maps give the
following gap below $\varphi(n)$.

\begin{proposition}\label{prop:gapCP}
Let $f$ be a congruence-preserving permutation of $\Z_n$. If
$|D_n(f)|<\varphi(n)$, then $|D_n(f)|\le G(n)$.
\end{proposition}

\begin{proof}
For every prime $p\mid n$, let $T_p=T(f_p)$ be as in the proof of
Theorem~\ref{thm:half}. Then, $|D_n(f)|=n\prod_{p\mid n}(|T_p|/p)$, and
$|T_p|\le p-1$ by Theorem~\ref{thm:universal}. If $|T_p|=p-1$ for every
$p$, then $|D_n(f)|=\varphi(n)$. Hence, there is a prime $\ell\mid n$ with
$|T_\ell|\le\ell-2$. Since $\ell\le q$, we obtain
$$
|D_n(f)|\le\varphi(n)\,\frac{\ell-2}{\ell-1}\le\varphi(n)\,\frac{q-2}{q-1}.
$$
\end{proof}

For permutations that are not congruence-preserving, we use the
following rigidity property of maps of $\Z_m\times\F_q$.

\begin{lemma}\label{lem:rowrigidity}
Let $m\ge2$ be an integer, let $q>m$ be a prime, and let
$g=(g_1,g_2):\Z_m\times\F_q\to\Z_m\times\F_q$ be a map.
Suppose that, for every $\beta\in\F_q^\times$, some $\alpha\in\Z_m$
makes
$$
(s,t)\longmapsto
\bigl(g_1(s,t)+\alpha s,\,g_2(s,t)+\beta t\bigr)
$$
a permutation. Fix $a\in\Z_m$. If this map is a permutation for
more than $(q-1)/2$ choices of $\beta\ne0$ with $\alpha=a$, then
$(s,t)\mapsto(g_1(s,t)+as,t)$ is a permutation, and there is a map
$h_a:\Z_m\to\F_q$ such that
$$
g_2(s,t)=h_a\bigl(g_1(s,t)+as\bigr).
$$
Consequently, with $\alpha=a$, the map is a permutation for every
$\beta\ne0$.
\end{lemma}

\begin{proof}
For each $\beta\ne0$, choose $\alpha$ as in the hypothesis. The second
coordinate $g_2(s,t)+\beta t$ of the resulting permutation takes every
value exactly $m$ times. Apply Lemma~\ref{lem:uniform}(ii) over
$\F_q$, with first variable $t$ and index set $\Z_m$. It follows that
the multiset $\{g_2(s,t):s\in\Z_m\}$ is independent of $t$.
In particular, the image of $g_2$ has at most $m$ elements.

Put $F_a(s,t)=g_1(s,t)+as$. By hypothesis, $(F_a,g_2+\beta t)$ is a
permutation for some $\beta$. Hence, each fibre $F_a^{-1}(u)$ has $q$
elements. The map $(s,t)\mapsto(t,g_2(s,t))$ is injective on this
fibre. Indeed, equal values would give a collision in this permutation.
Let $S_u\subseteq\F_q^2$ be its image. For every $\beta$ for which
$(F_a,g_2+\beta t)$ is a permutation, the map $(t,z)\mapsto z+\beta t$
is injective on $S_u$. In the terminology of Lemma~\ref{lem:pointsets},
the set $S_u$ is injective in the direction $\langle(1,-\beta)\rangle$. There
are more than $(q-1)/2$
such directions. By Lemma~\ref{lem:pointsets}(ii), applied to the prime
$q$, the set $S_u$ is a line.

For each fixed $t$, there are at most $m<q$ inputs $(s,t)$.
Hence, $S_u$ is not a vertical line. A non-horizontal, non-vertical
line would contain all $q$ values of its second coordinate.
This contradicts the bound $|\operatorname{im}(g_2)|\le m$.
Therefore, $S_u$ is a horizontal line. Its second coordinate has a
constant value $h_a(u)$, and its first coordinate runs over $\F_q$.
There is exactly one input in $F_a^{-1}(u)$ for each $t$.
This proves both assertions about $F_a$ and $h_a$.
On each fibre, $h_a(u)+\beta t$ is bijective for every $\beta\ne0$,
which proves the last assertion.
\end{proof}

We can now prove the gap bound for products of two distinct odd primes.

\begin{theorem}\label{thm:gappq}
Let $p<q$ be odd primes. If a permutation $f$ of $\Z_{pq}$ is not
extremal, then
$$
|D_{pq}(f)|\le(p-1)(q-2)=G(pq).
$$
\end{theorem}

\begin{proof}
Suppose that $|D_{pq}(f)|>(p-1)(q-2)$. We show that $f$ is extremal.
Choose a unit $c$ with $D_{pq}(f)\subseteq\Z_{pq}^\times-c$, as in
Theorem~\ref{thm:universal}, and put $g(x)=f(x)-cx$.
Let $U\subseteq\F_p^\times\times\F_q^\times$ consist of the pairs
$(\alpha,\beta)$ for which
$(s,t)\mapsto g(s,t)+(\alpha s,\beta t)$ is a permutation.
Then, $U=D_{pq}(f)+c$. Since $|U|>(p-1)(q-2)$, the number $K$ of pairs
of $\F_p^\times\times\F_q^\times$ outside $U$ satisfies $K\le p-2$.

Every nonzero $\alpha$ occurs in $U$, because $K<q-1$.
Every nonzero $\beta$ occurs as well, because $K<p-1$.
For $\alpha\in\F_p^\times$, we call the set of pairs $(\alpha,\beta)$ with
$\beta\in\F_q^\times$ the row of $\alpha$. Lemma~\ref{lem:rowrigidity},
applied with $m=p$, shows that either the whole row of $\alpha$ lies in
$U$, or at most $(q-1)/2$ of its pairs lie in $U$. Each incomplete row
therefore contains at least $(q-1)/2$ pairs outside $U$. Since $K<q-1$,
at most one row is incomplete. If all rows are complete, then
$|D_{pq}(f)|=\varphi(pq)$, and $f$ is extremal.

Suppose that the row of $a_0$ is incomplete.
For $p=3$, it would contain at least $(q-1)/2\ge2$ pairs outside $U$,
contrary to $K\le1$. Thus, $p\ge5$.
Put $S=\F_p^\times\setminus\{a_0\}$.
For $a\in S$, Lemma~\ref{lem:rowrigidity} gives a map $h_a$ with
$$
F_a(s,t)=g_1(s,t)+as,
\quad g_2(s,t)=h_a(F_a(s,t)),
$$
and $s\mapsto F_a(s,t)$ is a permutation of $\F_p$ for every fixed $t$.
Consequently, the map $s\mapsto g_1(s,t)$ has at least
$p-2>(p-3)/2$ good slopes over $\F_p$.
By Theorem~\ref{thm:fields}(i), it is affine. Write
$g_1(s,t)=A_t s+B_t$. Since $A_t+a\ne0$ for every $a\in S$,
we have $A_t\in\{0,-a_0\}$.

Suppose that $A_{t_1}=-a_0$ for some $t_1$.
There is at least one pair $(a_0,\beta_0)\in U$.
On the inputs $(s,t_1)$, the first coordinate
$F_{a_0}(s,t_1)=B_{t_1}$ of the corresponding permutation is constant. Hence,
$g_2(s,t_1)$ takes $p$ distinct values as $s$ varies.
For each $a\in S$, the map $s\mapsto F_a(s,t_1)$ is bijective.
Thus, $h_a$ is injective, and all $h_a$ have the same image, namely
the global image of $g_2$.

There must also be $t_0$ with $A_{t_0}=0$.
Otherwise, $F_{a_0}(s,t)=B_t$ for all $s,t$, and each fibre of
$F_{a_0}$ would have size divisible by $p$.
But $(a_0,\beta_0)\in U$ requires each fibre to have size $q$,
and $p\nmid q$.
Choose distinct $a,a'\in S$ and put $\pi=h_{a'}^{-1}\circ h_a$.
The inverses are taken on their common image.
Then, $F_{a'}(s,t)=\pi(F_a(s,t))$ for all $s,t$.
At $t=t_0$, this forces $\pi$ to be affine with multiplier $a'/a$.
At $t=t_1$, it forces its multiplier to be $(a'-a_0)/(a-a_0)$.
These cannot agree, since equality would imply
$a_0(a-a')=0$. This contradiction proves that $A_t=0$ for every $t$.

We now have $g_1(s,t)=B_t$. Fix distinct $a,a'\in S$.
For every $z\in\F_p$ and every $t\in\F_q$, the identities for $h_a$
and $h_{a'}$ give
$$
h_a(z)=h_{a'}\bigl((a'/a)z+(1-a'/a)B_t\bigr).
$$
If the values $B_t$ are not all equal, two choices of $t$ show that
$h_{a'}$ is invariant under a nonzero translation of $\F_p$.
Such a translation generates the additive group of $\F_p$. Hence, $h_{a'}$ is
constant. Then, $g_2$ is constant and $g_1$ depends only on $t$.
If all $B_t$ are equal, then $g_1$ is constant and
$g_2(s,t)=h_a(as+B_t)$ depends only on $s$.
In either situation, $g$ is strictly triangular with respect to a chain
$\Z_{pq}\supsetneq I\supsetneq\{0\}$, where $I$ is $\F_p\times\{0\}$ or
$\{0\}\times\F_q$. Proposition~\ref{prop:triangular}, applied to
$f(x)=cx+g(x)$, shows that $f$ is extremal. This completes the proof.
\end{proof}

Lemma~\ref{lem:rowrigidity} also settles the moduli whose largest prime
divisor is large compared with the rest of the modulus.

\begin{theorem}\label{thm:gapmq}
Let $n=mq$ be odd, where $m>1$ and $q$ is a prime with $q>m$ and
$\varphi(m)\le(q-1)/2$. Then, every non-extremal permutation $f$ of
$\Z_n$ satisfies $|D_n(f)|\le G(n)$.
\end{theorem}

\begin{proof}
Since $q>m$, the prime $q$ is the largest prime divisor of $n$, and
$\gcd(m,q)=1$. Hence, $\varphi(n)=\varphi(m)(q-1)$ and
$G(n)=\varphi(m)(q-2)$. Suppose that $|D_n(f)|>G(n)$. We show that $f$
is extremal. As in the proof of Theorem~\ref{thm:gappq}, identify $\Z_n$
with $\Z_m\times\F_q$, choose a unit $c$ with
$D_n(f)\subseteq\Z_n^\times-c$, and put $g(x)=f(x)-cx$. Let $U$ be the set
of units $(\alpha,\beta)$ for which $(s,t)\mapsto g(s,t)+(\alpha s,\beta t)$
is a permutation. Then, $U=D_n(f)+c$, and the number $K$ of units outside
$U$ satisfies $K<\varphi(m)(q-1)-\varphi(m)(q-2)=\varphi(m)$.

For every $\beta\in\F_q^\times$, there are $\varphi(m)>K$ units with second
coordinate $\beta$. Hence, every $\beta$ occurs in $U$. For every
$\alpha\in\Z_m^\times$, at most $K<\varphi(m)\le(q-1)/2$ values of $\beta$
satisfy $(\alpha,\beta)\notin U$. Thus, $(\alpha,\beta)\in U$ for more
than $(q-1)/2$ values of $\beta$, and Lemma~\ref{lem:rowrigidity} shows
that $(\alpha,\beta)\in U$ for every $\beta\in\F_q^\times$. Therefore,
$K=0$, $|D_n(f)|=\varphi(n)$, and $f$ is extremal.
\end{proof}

For instance, Theorem~\ref{thm:gapmq} applies to $n=9q$ for every prime
$q\ge13$, and to $n=15q$ for every prime $q\ge17$. For $m=p$ prime, it
is the case $q\ge2p-1$ of Theorem~\ref{thm:gappq}.

The gap bound also behaves well under the decomposition of
Lemma~\ref{lem:split-general}. The next estimate records the
contribution of the primes that do not divide the smaller modulus.

\begin{lemma}\label{lem:relativebound}
Let $f$ be a permutation of $\Z_n$, let $\ell\mid n$ with $\ell\ge2$,
and let $h$ be a permutation of $\Z_\ell$. If reduction modulo $\ell$
maps $D_n(f)$ into $D_\ell(h)$, then
$$
\frac{|D_n(f)|}{\varphi(n)}\le
\frac{|D_\ell(h)|}{\varphi(\ell)}.
$$
\end{lemma}

\begin{proof}
Choose a unit $a$ with $D_n(f)\subseteq\Z_n^\times-a$, using
Theorem~\ref{thm:universal}. Above a fixed residue modulo $\ell$,
there are $n/\ell$ elements of $\Z_n$. For each prime $p\mid n$ with
$p\nmid\ell$, membership in $D_n(f)$ requires $m\not\equiv-a\pmod p$.
By the Chinese remainder theorem, these restrictions leave exactly
$$
\frac{n}{\ell}\prod_{\substack{p\mid n\\p\nmid\ell}}
\left(1-\frac1p\right)=\frac{\varphi(n)}{\varphi(\ell)}
$$
possible elements above each residue. Ignoring any further restrictions
from primes dividing $\ell$ gives the asserted upper bound.
\end{proof}

The number $G(\ell)$ is defined by the same formula for every
$\ell\ge2$. For an odd prime $\ell$, it equals $\ell-2$. Hence, by
Theorem~\ref{thm:universal}, every permutation of $\Z_\ell$ is extremal
or has at most $G(\ell)$ good slopes. The next theorem reduces the gap
bound through a preserved congruence.

\begin{theorem}\label{thm:gapsplit}
Let $n$ be odd, let $n=dk$ with $d,k\ge2$, and let $f$ be a
permutation of $\Z_n$ that preserves congruences modulo $d$.
Let $g$ and $f_r$ be the permutations in
Lemma~\ref{lem:split-general}. Suppose that each of these permutations,
say $h$ on $\Z_\ell$ with $\ell\in\{d,k\}$, is extremal or satisfies
$|D_\ell(h)|\le G(\ell)$. Then, $f$ is extremal or satisfies
$|D_n(f)|\le G(n)$.
\end{theorem}

\begin{proof}
Let $q$ be the largest prime divisor of $n$. Suppose first that one
of the smaller permutations, denoted by $h$ on $\Z_\ell$, is not
extremal. Let $q_\ell$ be the largest prime divisor of $\ell$.
By Lemmas~\ref{lem:split-general} and~\ref{lem:relativebound},
$$
\frac{|D_n(f)|}{\varphi(n)}\le
\frac{|D_\ell(h)|}{\varphi(\ell)}
\le\frac{q_\ell-2}{q_\ell-1}\le\frac{q-2}{q-1}.
$$
This is the required gap bound.

Suppose now that $g$ and all $f_r$ are extremal. By
Theorem~\ref{thm:universal}, the set of good slopes of each of them is
$\Z_\ell^\times-a$ for some unit $a$ of the ring $\Z_\ell$ on which it
acts. For every prime $p\mid n$, let $A_p\subseteq\F_p^\times$ be the
set of the residues modulo $p$ of these units $a$, taken over the
permutations whose modulus $\ell$ is divisible by $p$. Each $A_p$ is
nonempty, because $p$ divides $d$ or $k$. Since $\Z_\ell^\times-a$
consists of the $m$ with $m\not\equiv-a\pmod p$ for every prime
$p\mid\ell$, Lemma~\ref{lem:split-general} gives
$$
D_n(f)=\{m\in\Z_n:m\bmod p\notin-A_p\text{ for every }p\mid n\},
\quad |D_n(f)|=n\prod_{p\mid n}\frac{p-|A_p|}{p}.
$$
If every $A_p$ has one element, this cardinality is $\varphi(n)$.
Otherwise, $|A_p|\ge2$ for some prime $p$, and
$$
|D_n(f)|\le\varphi(n)\frac{p-2}{p-1}
\le\varphi(n)\frac{q-2}{q-1}=G(n).
$$
This proves the theorem.
\end{proof}

\begin{corollary}\label{cor:minimalgap}
Let $n$ be odd and composite, and suppose that some non-extremal
permutation of $\Z_n$ has more than $G(n)$ good slopes. Then, there are
a composite divisor $\ell$ of $n$ and a non-extremal permutation $h$ of
$\Z_\ell$ with $|D_\ell(h)|>G(\ell)$ that preserves congruences modulo no
$d$ with $1<d<\ell$.
\end{corollary}

\begin{proof}
We call a non-extremal permutation $h$ of $\Z_\ell$ a counterexample if
$|D_\ell(h)|>G(\ell)$. Start with a counterexample on $\Z_n$. If it preserves
congruences modulo some $d$ with $1<d<n$, then
Theorem~\ref{thm:gapsplit} shows that one of its quotient or
residue-class permutations is a counterexample. Its modulus is a proper
divisor of $n$. We repeat this step as long as a nontrivial congruence
is preserved. The modulus decreases strictly, and the process stops. By
Theorem~\ref{thm:universal}, there is no counterexample of prime order.
Hence, the final permutation has the stated properties.
\end{proof}

\begin{proposition}\label{prop:gapsharp}
Let $n$ be odd and composite, and let $q$ be its largest prime divisor.
Then, there is a non-affine permutation $f$ of $\Z_n$ with
$|D_n(f)|=G(n)$. If $q^2\mid n$, such a permutation can be chosen to be
polynomial.
\end{proposition}

\begin{proof}
If $q^2\mid n$, the polynomial $F(X)=X+c(X-X^q)$ of Example~\ref{ex:2xx5}
has the required properties. Now suppose that $q\,\|\,n$, and write $n=mq$
with $m>1$ and $\gcd(m,q)=1$. Identify $\Z_n$ with $\Z_m\times\F_q$. Choose
a map $\lambda:\Z_m\to\F_q^\times$ with exactly two distinct values
$\lambda_1$ and $\lambda_2$, which is possible because $m\ge2$ and
$q\ge3$. Put $f(s,t)=(s,\lambda(s)t)$. For a slope
$(\alpha,\beta)\in\Z_m\times\F_q$, we have
$$
f(s,t)+(\alpha,\beta)(s,t)=\bigl((1+\alpha)s,\,(\lambda(s)+\beta)t\bigr).
$$
If $1+\alpha$ is not a unit of $\Z_m$, then the first coordinate does not
take all values. If $1+\alpha$ is a unit, then the first coordinate
determines $s$, and the map is a permutation exactly when
$\lambda(s)+\beta\ne0$ for all $s$, that is, when
$\beta\notin\{-\lambda_1,-\lambda_2\}$. In particular, $f$ is a
permutation, and $|D_n(f)|=\varphi(m)(q-2)=G(n)$. In both constructions,
$|D_n(f)|<\varphi(n)$, while every affine permutation has $\varphi(n)$ good
slopes. Hence, the permutations are not affine.
\end{proof}

\begin{example}\label{ex:gap15}
On $\Z_{15}=\F_3\times\F_5$, take $f(s,t)=\bigl(s,(1+[s=0])t\bigr)$.
This is the permutation $f(s,t)=(s,\lambda(s)t)$ of the proof of
Proposition~\ref{prop:gapsharp} with $\lambda(0)=2$ and $\lambda(s)=1$ for
$s\neq0$. Its good slopes $(\alpha,\beta)$ are those with
$\alpha\in\{0,1\}$ and $\beta\in\{0,1,2\}$. Hence,
$$
D_{15}(f)=\{0,1,6,7,10,12\},\quad
|D_{15}(f)|=6=\varphi(15)\,\frac{5-2}{5-1},
$$
and $\mu_{15}(f)=2$. This permutation is not congruence-preserving. Indeed,
the elements $(0,1)$ and $(1,1)$ are congruent modulo $5$, but their images
$(0,2)$ and $(1,1)$ are not.
\end{example}

Figure~\ref{fig:crt-slopes-fifteen} compares this slope set with the
extremal one in Example~\ref{ex:coprime}. By
Lemma~\ref{lem:admissible}(i), an admissible subset corresponds to cells
in distinct rows and distinct columns. In both panels, the good slopes
occupy only two rows. Hence, an admissible subset has at most two
elements, and the pair $\{0,1\}$ attains this bound.

\begin{figure}[htbp]
\centering
\begin{tikzpicture}[x=.94cm,y=.75cm,every node/.style={font=\small}]
\begin{scope}[shift={(0,0)}]
\fill[black!12] (-0.5,1.5) rectangle (0.5,2.5);
\fill[black!12] (0.5,1.5) rectangle (1.5,2.5);
\fill[black!12] (1.5,1.5) rectangle (2.5,2.5);
\fill[black!12] (2.5,1.5) rectangle (3.5,2.5);
\fill[black!12] (-0.5,0.5) rectangle (0.5,1.5);
\fill[black!12] (0.5,0.5) rectangle (1.5,1.5);
\fill[black!12] (1.5,0.5) rectangle (2.5,1.5);
\fill[black!12] (2.5,0.5) rectangle (3.5,1.5);
\foreach \i in {0,...,5}{
  \draw[black!45,line width=.35pt] ({\i-.5},-.5) -- ({\i-.5},2.5);
}
\foreach \j in {0,...,3}{
  \draw[black!45,line width=.35pt] (-.5,{\j-.5}) -- (4.5,{\j-.5});
}
\draw[black!65,line width=.5pt] (-.5,-.5) rectangle (4.5,2.5);
\foreach \b in {0,...,4}{\node at (\b,3.05) {$\b$};}
\foreach \a in {0,1,2}{\node at (-1,{2-\a}) {$\a$};}
\node at (2,3.8) {$\beta=m\bmod5$};
\node[rotate=90] at (-1.8,1) {$\alpha=m\bmod3$};
\node at (2,4.65) {(a) Example~\ref{ex:coprime}};
\node at (0,2) {$0$};
\draw[black,line width=.8pt] (0,2) circle[radius=.23cm];
\node at (1,2) {$6$};
\node at (2,2) {$12$};
\node at (3,2) {$3$};
\node at (4,2) {$9$};
\node at (0,1) {$10$};
\node at (1,1) {$1$};
\draw[black,line width=.8pt] (1,1) circle[radius=.23cm];
\node at (2,1) {$7$};
\node at (3,1) {$13$};
\node at (4,1) {$4$};
\node at (0,0) {$5$};
\node at (1,0) {$11$};
\node at (2,0) {$2$};
\node at (3,0) {$8$};
\node at (4,0) {$14$};
\node at (2,-1.2) {$|D_{15}(f)|=8,\quad\mu_{15}(f)=2$};
\end{scope}
\begin{scope}[shift={(8,0)}]
\fill[black!12] (-0.5,1.5) rectangle (0.5,2.5);
\fill[black!12] (0.5,1.5) rectangle (1.5,2.5);
\fill[black!12] (1.5,1.5) rectangle (2.5,2.5);
\fill[black!12] (-0.5,0.5) rectangle (0.5,1.5);
\fill[black!12] (0.5,0.5) rectangle (1.5,1.5);
\fill[black!12] (1.5,0.5) rectangle (2.5,1.5);
\foreach \i in {0,...,5}{
  \draw[black!45,line width=.35pt] ({\i-.5},-.5) -- ({\i-.5},2.5);
}
\foreach \j in {0,...,3}{
  \draw[black!45,line width=.35pt] (-.5,{\j-.5}) -- (4.5,{\j-.5});
}
\draw[black!65,line width=.5pt] (-.5,-.5) rectangle (4.5,2.5);
\foreach \b in {0,...,4}{\node at (\b,3.05) {$\b$};}
\foreach \a in {0,1,2}{\node at (-1,{2-\a}) {$\a$};}
\node at (2,3.8) {$\beta=m\bmod5$};
\node[rotate=90] at (-1.8,1) {$\alpha=m\bmod3$};
\node at (2,4.65) {(b) Example~\ref{ex:gap15}};
\node at (0,2) {$0$};
\draw[black,line width=.8pt] (0,2) circle[radius=.23cm];
\node at (1,2) {$6$};
\node at (2,2) {$12$};
\node at (3,2) {$3$};
\node at (4,2) {$9$};
\node at (0,1) {$10$};
\node at (1,1) {$1$};
\draw[black,line width=.8pt] (1,1) circle[radius=.23cm];
\node at (2,1) {$7$};
\node at (3,1) {$13$};
\node at (4,1) {$4$};
\node at (0,0) {$5$};
\node at (1,0) {$11$};
\node at (2,0) {$2$};
\node at (3,0) {$8$};
\node at (4,0) {$14$};
\node at (2,-1.2) {$|D_{15}(f)|=6,\quad\mu_{15}(f)=2$};
\end{scope}

\end{tikzpicture}
\caption{Good slopes modulo $15$, arranged by their residues modulo
$3$ and $5$. Each cell contains the unique residue $m\in\Z_{15}$
with the indicated pair of residues. Shading marks good slopes, and
circles mark the admissible set $\{0,1\}$. Panel (a) attains
$\varphi(15)=8$. Panel (b) attains the gap bound $G(15)=6$.}
\label{fig:crt-slopes-fifteen}
\end{figure}

For $n=15$, we have $G(15)=6<15/2$. Thus, by Theorem~\ref{thm:gappq},
every permutation of $\Z_{15}$ with more than $15/2$ good slopes is
extremal, and
Theorem~\ref{thm:pq} describes it. The same conclusion holds for
$n=21$, because $G(21)=10<21/2$.

\section{The R\'edei polynomial over \texorpdfstring{$\Z_n$}{Zn}}
\label{sec:redei}

R\'edei's method studies the directions determined by a function over a
finite field through a polynomial in two variables. In this section, we
consider the analogous polynomial over $\Z_n$. It gives a second proof of
the first assertion of Theorem~\ref{thm:universal}. We also explain why
the second half of R\'edei's method has no direct analogue over $\Z_n$.

For a permutation $f$ of $\Z_n$, consider the polynomials
$$
F_f(X,Y)=\prod_{x\in\Z_n}\bigl(X+xY+f(x)\bigr)\in\Z_n[X,Y]
\quad\text{and}\quad
P_n(X)=\prod_{\alpha\in\Z_n}(X+\alpha)\in\Z_n[X].
$$
They are polynomials over $\Z_n$ even when $f$ is not a polynomial map. If
$m\in D_n(f)$, then the factors $X+f(x)+mx$ are the factors of $P_n$ in
another order. Hence, $F_f(X,m)=P_n(X)$ for every $m\in D_n(f)$. The same
definitions make sense for a map $f:\F_q\to\F_q$. There, the polynomial
$P_q(X)=X^q-X$ has $q$ distinct linear factors, and unique factorisation in
$\F_q[X]$ shows that $F_f(X,m)=P_q(X)$ holds if and only if
$m\in D_{\F_q}(f)$. This equivalence is the starting point of R\'edei's
method. Write
$$
F_f(X,Y)=X^n+\sum_{j=1}^nC_j(Y)X^{n-j}\quad\text{and}\quad
P_n(X)=X^n+\sum_{j=1}^nb_jX^{n-j}.
$$
Here, $C_j(Y)$ is the $j$th elementary symmetric polynomial of the $n$
linear forms $f(x)+xY$, and $b_j$ is the $j$th elementary symmetric
polynomial of $0,1,\dots,n-1$. Hence, $\deg C_j\le j$, the coefficient of
$Y^j$ in $C_j$ is $b_j$, and $C_j(m)=b_j$ for every $m\in D_n(f)$.

\begin{proposition}\label{prop:redei}
Let $f$ be a permutation of $\Z_n$.
\begin{enumerate}
\item If $M\subseteq D_n(f)$ is admissible and $j<|M|$, then $C_j(Y)=b_j$.
\item Let $p$ be a prime with $p^e\,\|\,n$. If $p$ is odd, then
$b_j\equiv0\pmod{p^e}$ for $1\le j\le p-2$, and
$b_{p-1}\equiv-n/p\pmod{p^e}$. If $p=2$, then $b_1\equiv-n/2\pmod{2^e}$.
\item In the situation of (ii), the polynomial $C_{p-1}(Y)-b_{p-1}$
vanishes on $D_n(f)$, and modulo $p^e$ its coefficient of $Y^{p-1}$ is
$b_{p-1}\not\equiv0$. Consequently, $D_n(f)$ does not meet every residue
class modulo $p$.
\end{enumerate}
\end{proposition}

\begin{proof}
For (i), the polynomial $C_j(Y)-b_j$ has degree at most $j<|M|$ and
vanishes on $M$. By Lemma~\ref{lem:interp}, it is zero.

For (ii), let $s_i=\sum_{x=0}^{n-1}x^i$. Newton's identities for the integers
$0,1,\dots,n-1$ read $kb_k=\sum_{i=1}^k(-1)^{i-1}b_{k-i}s_i$, where $k\ge1$
and $b_0=1$. For $p=2$, this gives $b_1=s_1\equiv-n/2\pmod{2^e}$
by Lemma~\ref{lem:powersum}. Let $p$ be odd. By Lemma~\ref{lem:powersum},
$s_i\equiv0\pmod{p^e}$ for $1\le i\le p-2$. For $1\le k\le p-2$, every term
on the right contains such an $s_i$, and $k$ is a unit modulo $p^e$. Hence,
$b_k\equiv0\pmod{p^e}$. For $k=p-1$, the only term that
does not contain one of $s_1,\dots,s_{p-2}$ is $(-1)^{p-2}s_{p-1}=-s_{p-1}$.
Thus, $(p-1)b_{p-1}\equiv-s_{p-1}\equiv n/p\pmod{p^e}$. Consequently,
$(p-1)(b_{p-1}+n/p)\equiv n/p+(p-1)n/p=n\equiv0\pmod{p^e}$. Since $p-1$ is
a unit modulo $p^e$, we obtain $b_{p-1}\equiv-n/p\pmod{p^e}$.

For (iii), the polynomial $C_{p-1}(Y)-b_{p-1}$ vanishes on $D_n(f)$,
because $C_{p-1}(m)=b_{p-1}$ for $m\in D_n(f)$. The assertion about its
leading coefficient follows from (ii), because $n/p$ is not divisible by
$p^e$. If $D_n(f)$ met every residue class modulo $p$, it would contain
$p$ elements whose differences are units modulo $p^e$. Then,
Lemma~\ref{lem:interp}, applied over $\Z_{p^e}$, would force
$b_{p-1}\equiv0\pmod{p^e}$, a contradiction.
\end{proof}

Proposition~\ref{prop:redei}(iii) gives a second proof of the first
assertion of Theorem~\ref{thm:universal}, directly from the R\'edei
polynomial. It also explains why we work modulo the full prime power
$p^e$. Modulo $p$ alone, the leading coefficient $-n/p$ vanishes as soon
as $p^2\mid n$. The second half of R\'edei's method uses unique
factorisation of $F_f(X,y)$ when $y$ is a determined direction. Over
$\Z_n$, the identity $F_f(X,m)=P_n(X)$ does not characterise good slopes,
as the next example shows.

\begin{example}\label{ex:product}
In $\Z_9[X]$, we have
$(X+r)(X+r+3)(X+r+6)=(X+r)^3+9(X+r)^2+18(X+r)=(X+r)^3$. We group the
elements of $\Z_9$ according to their residues $0$, $1$ and $-1$ modulo
$3$. This gives $P_9(X)=\bigl(X(X+1)(X-1)\bigr)^3=(X^3-X)^3$. The map
$x\mapsto x^3$ of $\Z_9$ takes each of the values $0$, $1$ and
$-1$ exactly three times. Hence, $\prod_{x\in\Z_9}(X+x^3)=P_9(X)$, although
$x\mapsto x^3$ is not a permutation. The same phenomenon occurs for a
permutation and a slope that is not good. Let $f(x)=x^3-2x$. Then,
$f(x)-f(y)=(x-y)(x^2+xy+y^2-2)$, and modulo $3$ the second factor equals
$(x-y)^2+1$, which never vanishes. Hence, the second factor is a unit of
$\Z_9$, and $f$ is a permutation. Proposition~\ref{prop:poly} gives
$D_9(f)=\{0,3,6\}$, because $\overline F_3(x)=-x$ and $F'\equiv1\pmod3$.
Nevertheless, $f(x)+2x=x^3$, and therefore $F_f(X,2)=P_9(X)$ although
$2\notin D_9(f)$.

For squarefree $n$, the situation is easy to describe. Then,
$\Z_n[X]\cong\prod_{p\mid n}\F_p[X]$ and $P_n\equiv(X^p-X)^{n/p}\pmod p$.
Hence, $F_f(X,m)=P_n(X)$ holds if and only if, for every $p\mid n$, each
residue class modulo $p$ contains exactly $n/p$ of the values $f(x)+mx$.
For $n=6$, for instance, the multiset $\{0,0,1,1,2,5\}$ has the same
product $\prod(X+\alpha)$ as $\Z_6$.
\end{example}

\section{Application to quantum Latin squares}\label{sec:qls}

In this section, we apply our results to the construction of mutually
orthogonal quantum Latin squares of Ball and Simoens \cite{BS}. We
determine the largest number of squares that this construction gives over
$\Z_n$, and we compare it with the construction over other finite rings.
We first recall the definitions from \cite{BS}. A \emph{quantum Latin square} of
order $n$ is an $n\times n$ array $\Psi=(\psi_{ij})$ of unit vectors of
$\C^n$ such that the entries of every row and of every column form an
orthonormal basis of $\C^n$. Two quantum Latin squares $\Psi$ and $\Phi$
are \emph{orthogonal} if the $n^2$ vectors $\psi_{ij}\otimes\phi_{ij}$ form
an orthonormal basis of $\C^n\otimes\C^n$. A set of $t$ pairwise
orthogonal quantum Latin squares of order $n$ is called a set of $t$
MOQLS$(n)$. A quantum Latin square is \emph{classical} if there is a
unitary map $U$ of $\C^n$ such that every vector $U\psi_{ij}$ is a scalar
multiple of a vector of the standard basis.

In \cite[Theorem~6]{BS}, the construction is stated for a finite
Frobenius ring $R$, and $\Z_n$ is such a ring. It uses a map $f$ that is not
the sum of an additive map and a constant, together with a set $M$ of
elements with pairwise invertible differences such that
$x\mapsto f(x)-mx$ is a permutation for every $m\in M$. Replacing $f$
by $x\mapsto f(x)-m_0x$ and $M$ by $M-m_0$ for some $m_0\in M$, we may
assume that $f$ is a permutation. The additive maps of $\Z_n$ are the
maps $x\mapsto ax$, and $x\mapsto f(x)-mx$ is a permutation exactly when
$-m\in D_n(f)$. Hence, $-M$ is an admissible subset of $D_n(f)$, and, for
$R=\Z_n$, the theorem reads as follows in our notation.

\begin{theorem}[{\cite[Theorem~6]{BS}}]\label{thm:qls}
Let $f$ be a non-affine permutation of $\Z_n$, and let
$M\subseteq D_n(f)$ be a nonempty admissible set. Then, there exist
$|M|$ MOQLS$(n)$, exactly one of which is non-classical.
\end{theorem}

The squares are explicit. Let $\omega=e^{2\pi i/n}$, and index the
standard basis of $\C^n$ by $\Z_n$. By Lemma~\ref{lem:normal}, we may
replace $f$ by $f(x)+m_0x$ and $M$ by $M-m_0$ for some $m_0\in M$. Thus, we
may assume that $0\in M$, and then every nonzero element of $M$ is a unit.
The non-classical square $\Psi$ has the entries
$$
\psi_{ij}=\frac1{\sqrt n}\sum_{x\in\Z_n}\omega^{if(x)+jx}\ket x,
$$
where $i,j\in\Z_n$. For every $m\in M\setminus\{0\}$, the classical square
$\Phi^{(m)}$ has the entries $\ket{j-mi}$. The square $\Psi$ is orthogonal
to $\Phi^{(m)}$ because $x\mapsto f(x)+mx$ is a permutation. The squares
$\Phi^{(m)}$ and $\Phi^{(m')}$ are orthogonal because $m-m'$ is a unit.
The square $\Psi$ is non-classical because $f$ is not affine.

Together with Theorem~\ref{thm:maxima}, Theorem~\ref{thm:qls} gives the
largest number of squares that the construction provides over $\Z_n$.

\begin{corollary}\label{cor:qlsZn}
Let $n\ge4$. The largest number of MOQLS$(n)$ that
Theorem~\ref{thm:qls} provides over $\Z_n$ is $B(n)$. It equals
$p_0(n)-1$ for composite $n\ge6$, it equals $(p-3)/2$ for a prime $p\ge5$,
and it equals $1$ for $n=4$. In particular, for even $n$, the construction
over $\Z_n$ provides a single square.
\end{corollary}

\begin{proof}
The number of squares in Theorem~\ref{thm:qls} is $|M|$. By definition,
the largest size of an admissible subset of $D_n(f)$, over all non-affine
permutations $f$ of $\Z_n$, is $B(n)$. Theorem~\ref{thm:maxima} gives its
value. For even $n\ge6$, we have $p_0(n)-1=1$.
\end{proof}

\begin{example}\label{ex:qls25}
Let $n=25$ and $f(x)=x+5x^2$. By Example~\ref{ex:primepower}, $D_{25}(f)$
consists of the elements $m\in\Z_{25}$ with $m\not\equiv4\pmod5$, and
$M=\{0,1,2,3\}$ is
admissible. With $\omega=e^{2\pi i/25}$, Theorem~\ref{thm:qls} gives the
non-classical square with entries
$$
\psi_{ij}=\frac15\sum_{x\in\Z_{25}}\omega^{\,i(x+5x^2)+jx}\ket x,
$$
together with the three classical squares with entries $\ket{j-mi}$ for
$m=1,2,3$. These are four MOQLS$(25)$, exactly one of which is
non-classical. By Corollary~\ref{cor:qlsZn}, no non-affine permutation of
$\Z_{25}$ gives more. Over $\F_{25}$, the same construction gives eleven
such squares \cite[Corollary~7]{BS}.
\end{example}

We finally prove a bound for finite rings whose order is not a prime
power. It applies, in particular, to every such ring allowed in the
construction of \cite{BS}. Here, $R$ is any finite ring with identity,
not necessarily commutative, and $R^\times$ is its group of units.

\begin{proposition}\label{prop:rings}
Let $R$ be a finite ring with identity whose order
$n=p_1^{a_1}\cdots p_k^{a_k}$ is not a prime power. Let $f:R\to R$ be a
permutation, and let $M\subseteq R$ be a set whose pairwise differences lie in
$R^\times$, such that $x\mapsto f(x)+mx$ is a permutation of $R$ for every
$m\in M$. Then, $|M|\le\min_{1\le i\le k}\bigl(p_i^{a_i}-1\bigr)$.
\end{proposition}

\begin{proof}
Let $R_i$ be the set of elements of $R$ whose additive order is a power of
$p_i$. Each $R_i$ is a two-sided ideal, and $R_iR_j=0$ for $i\neq j$.
Writing $1=\sum_ie_i$ with $e_i\in R_i$, we obtain a ring isomorphism
$R\cong\prod_iR_i$, where $R_i$ has identity $e_i$ and order $p_i^{a_i}$.
An element of $R$ is a unit if and only if all its components are units.
Fix $i$, write $S=R_i$, $q=p_i^{a_i}$ and $p=p_i$, and let $\pi:R\to S$ be
the projection. For distinct $m,m'\in M$, the element
$\pi(m)-\pi(m')=\pi(m-m')$ is a unit of $S$. In particular, $\pi$ is
injective on $M$.

Suppose first that $S$ has a nonzero element $s$ that is not a unit. We
claim that the map $x\mapsto xs$ of $S$ is not surjective. Otherwise,
$xs=1$ for some $x$. Then, $y\mapsto sy$ is injective, because $sy=sy'$
implies $y=xsy=xsy'=y'$. Hence, $sy=1$ for some $y$, and
$x=x(sy)=(xs)y=y$ is an inverse of $s$, a contradiction. As $S$ is
finite, the map $x\mapsto xs$ is not injective either. Thus,
$K=\{x\in S:xs=0\}$ is a nonzero additive subgroup of $S$, and $|K|\ge p$.
If two elements of $\pi(M)$ lay in the same coset of $K$, their difference
would be a unit $w$ with $ws=0$, and this would force $s=0$. Therefore,
$|M|=|\pi(M)|\le|S|/|K|\le q/p\le q-1$.

Otherwise, every nonzero element of $S$ is a unit, and by Wedderburn's
theorem $S$ is the field $\F_q$. Since $n$ is not a prime power, $L=n/q>1$
and $p\nmid L$. Let $m\in M$. The map $x\mapsto f(x)+mx$ is a permutation
of $R$, and $\pi$ is $L$-to-one. Hence, the map
$x\mapsto\pi(f(x)+mx)=\pi(f(x))+\pi(m)\pi(x)$ takes every value of $S$
exactly $L$ times. Since $\sum_{\sigma\in S}\sigma^{q-1}=q-1=-1$ in $S$,
the polynomial
$$
Q(T)=\sum_{x\in R}\bigl(\pi(f(x))+T\pi(x)\bigr)^{q-1}+L\in S[T]
$$
vanishes at $\pi(m)$ for every $m\in M$. Its coefficient of $T^{q-1}$ is
$\sum_{x\in R}\pi(x)^{q-1}=-L$, which is nonzero in $S$. If $|M|\ge q$,
then $\pi(M)=S$, and the nonzero polynomial $Q$ of degree $q-1$ would have
$q$ roots. Hence, $|M|\le q-1$ in this case as well.
\end{proof}

MacNeish \cite{MacNeish} constructed $\min_i p_i^{a_i}-1$ mutually
orthogonal Latin squares of order $n$, and every Latin square gives a
classical quantum Latin square. Proposition~\ref{prop:rings} shows that,
for orders that are not prime powers, the construction of
\cite[Theorem~6]{BS} never gives more squares than this number, for every
finite ring covered by that construction. Indeed, that construction
uses a set $M$ with pairwise invertible differences such that
$x\mapsto f(x)-mx$ is a permutation for every $m\in M$. Choose
$m_0\in M$ and put $g(x)=f(x)-m_0x$. Then, $g$ is a permutation, the
set $m_0-M$ has pairwise invertible differences, and
$g(x)+(m_0-m)x=f(x)-mx$ for every $m\in M$. Thus,
Proposition~\ref{prop:rings} applies to $g$ and $m_0-M$.

For an odd prime $p$ and $e\ge2$, put $q=p^e$.
By Theorem~\ref{thm:fields}(i) and Proposition~\ref{prop:megyesi}, the
largest number supplied by the construction over $\F_q$ is $(q-3)/2$.
This is the number in \cite[Corollary~7]{BS}.
It exceeds $p-1$, the largest number supplied over $\Z_{p^e}$.
All these upper bounds concern this particular construction.
They are not upper bounds on arbitrary sets of quantum Latin squares.

\section{Further questions}\label{sec:questions}

Theorems~\ref{thm:classification} and~\ref{thm:pq} classify the
extremal permutations for prime powers and products of two distinct
odd primes. Theorem~\ref{thm:nontriangular} shows that this
description through chains of ideals fails for every other modulus.

\begin{question}\label{q:classify}
Classify the extremal permutations of $\Z_n$ when $n$ is neither a
prime power nor a product of two distinct odd primes.
\end{question}

Lemma~\ref{lem:gcd} gives an exact condition in terms of the
greatest common divisors of pairwise differences. The problem is to
turn this condition into a structural description and a counting
formula. The examples in Section~\ref{sec:mixed} show why a single
chain of ideals is insufficient.

\begin{question}\label{q:gap}
Let $n$ be odd with at least two distinct prime divisors, and let
$q$ be its largest prime divisor. Does every non-extremal permutation
$f$ of $\Z_n$ satisfy $|D_n(f)|\le\varphi(n)(q-2)/(q-1)$?
\end{question}

Proposition~\ref{prop:gapCP} proves this bound for
congruence-preserving permutations, and Theorem~\ref{thm:gappq}
proves it for products of two distinct odd primes.
Theorem~\ref{thm:gapmq} also covers products $mq$ under the
size conditions stated there.
Theorem~\ref{thm:gapsplit} gives a conditional reduction to smaller
moduli. By Corollary~\ref{cor:minimalgap}, if the bound fails for $n$,
then it fails for a permutation of $\Z_\ell$ that preserves
congruences modulo no $d$ with $1<d<\ell$, where $\ell$ is a composite
divisor of $n$.
Proposition~\ref{prop:gapsharp} proves that the bound, if valid,
is sharp for every odd composite modulus. In the companion paper, we
prove the bound for odd prime powers \cite[Corollary~9.7]{Gaps}. This
explains the restriction in Question~\ref{q:gap} to moduli with at
least two distinct prime divisors. There, we also prove the bound for
$n=9q$ with a prime $q\ge5$ and for $n=pq^e$ with odd primes $p<q$
such that $q\ge2p-1$ \cite[Theorems~9.13 and~9.16]{Gaps}. Together
with the results of this paper, they settle every odd composite
modulus below $105=3\cdot5\cdot7$.

The restriction to odd $n$ is necessary. For $n=12$, we have
$G(12)=2$. However, the permutation $f$ of $\Z_{12}$ defined by
$$
\bigl(f(0),f(1),\dots,f(11)\bigr)=(0,1,2,3,5,4,7,9,6,11,8,10)
$$
has $D_{12}(f)=\{0,2,8\}$. It is not extremal, since $\varphi(12)=4$.

\section*{Acknowledgements}
The author is grateful to Simeon Ball for suggesting
the problems that motivated this work.

{\small
\bibliographystyle{abbrv}
\bibliography{ref}
}

\end{document}